\documentclass[a4paper,11pt]{article}
\usepackage{latexsym,amsfonts,amsmath}
\usepackage{amssymb}
\usepackage{color,soul} 
\usepackage{xcolor}
\usepackage[mathscr]{eucal}
\usepackage{enumerate}
\usepackage{enumitem}

\usepackage{pgf}
\usepackage{tikz}
\usetikzlibrary{arrows,automata}

\usepackage{fullpage}

\definecolor{Dgreen}{rgb}{102,255,102}
\definecolor{Purp}{rgb}{102,0,204}

\def\1{\mathbf{1}}
\def\0{\mathbf{0}}

\def\NN{\mathbb{N}}

\def\RR{\mathbb{R}}

\def\XX{\mathbf{X}}
\def\AA{\mathbf{A}}
\def\KK{\mathbf{K}}
\def\HH{\mathbf{H}}

\def\MM{\boldsymbol{\mathcal{M}}}

\newcommand{\bcal}[1]{\boldsymbol{\mathcal{#1}}}

\newcommand{\bfrak}[1]{\boldsymbol{\mathfrak{#1}}}

\definecolor{Dgreen}{rgb}{0,0.5,0}

\newcounter{hypot}

\newcounter{assump}
\renewcommand{\theassump}{$(\mathrm{\Alph{hypot}}_\arabic{assump})$}

\newcounter{condS}
\renewcommand{\thecondS}{$(\mathrm{S}_\arabic{condS})$}
\newenvironment{conditS}{\begin{list}
      {\hspace{\labelsep} \bfseries \thecondS}
      {\leftmargin=0pt
       \labelwidth=0cm
       \usecounter{condS}
       \def\makelabel##1{##1}
       }}{\end{list}}
       
\newcounter{condW}
\renewcommand{\thecondW}{$(\mathrm{W}_\arabic{condW})$}
\newenvironment{conditW}{\begin{list}
      {\hspace{\labelsep} \bfseries \thecondW}
      {\leftmargin=0pt
       \labelwidth=0cm
       \usecounter{condW}
       \def\makelabel##1{##1}
       }}{\end{list}}       
       
\begin{document}
\newtheorem{theorem}{Theorem}[section]
\newtheorem{proposition}[theorem]{Proposition}
\newtheorem{condition}[theorem]{Condition}
\newtheorem{lemma}[theorem]{Lemma}
\newtheorem{corollary}[theorem]{Corollary}
\newtheorem{definition}[theorem]{Definition}
\newtheorem{remark}[theorem]{Remark}
\newtheorem{example}[theorem]{Example}
\newtheorem{conjecture}[theorem]{Conjecture}
\newtheorem{assumption}[theorem]{Assumption}

\bibliographystyle{plain}

\title{A note on weak compactness of occupation measures \\ for an absorbing Markov decision process
 \footnote{Supported by grant PID2021-122442NB-I00 from the Spanish \textit{Ministerio de Ciencia e Innovaci\'on.}}}

\author{Fran\c{c}ois Dufour\footnote{
Institut Polytechnique de Bordeaux; INRIA Bordeaux Sud Ouest, Team: ASTRAL; IMB, Institut de Math\'ematiques de Bordeaux, Universit\'e de Bordeaux, France \tt{francois.dufour@math.u-bordeaux.fr}}
\and
Tom\'as Prieto-Rumeau\footnote{Statistics Department, UNED, Madrid, Spain. e-mail: {\tt{tprieto@ccia.uned.es}}\quad {(Author for correspondence)}}}
\maketitle
\begin{abstract}  We consider an absorbing Markov decision process
with Borel state and action spaces. We study conditions under which the  MDP is uniformly absorbing and the set of occupation measures of the MDP is compact in the usual weak topology. These   include suitable continuity requirements on the transition kernel and conditions on the dynamics of the system at the boundary of the absorbing set. We generalize previously known results  and give an answer to some conjectures that have been mentioned in the related literature.
\end{abstract}
{\small 
\par\noindent\textbf{Keywords:} Markov decision processes; absorbing model; compactness of occupation measures
\par\noindent\textbf{AMS 2020 Subject Classification:} 90C40}

\section{Introduction}

We are interested in absorbing discrete-time Markov decision processes (MDPs) with Borel state and action spaces $\XX$ and $\AA$. An MDP is called absorbing when the controlled process reaches some absorbing subset $\Delta$ of the state space in a finite expected time for any policy. 
When considering this MDP under a total expected cost  optimality criterion, one of the most widely used techniques is the so-called linear programming or convex analytic approach: this consists in formulating the control problem as the equivalent linear optimization problem of minimizing a  linear cost functional over the set of \emph{occupation measures}. Roughly speaking, the occupation measure of a policy is a measure on $\XX\times\AA$ which computes, for each $\Gamma\subseteq\XX\times\AA$, the expected number of times that the state-action process visits $\Gamma$ before absorption by $\Delta$ takes place. 
Therefore, when using this approach, a key result is to prove that the set of occupation measures is compact for a suitable topology, under which the functional to be minimized is lower semicontinuous. This issue has been studied in, e.g., \cite{DuPri24-ESAIM,Feinberg-Piunovskiy-2019,feinberg12, piunovskiy2023continuity,schal75}, among others.

In general, the property of the MDP being absorbing does not suffice to ensure compactness of the set of occupation measures, and further hypotheses become necessary. In this context, the notion of a  \emph{uniformly absorbing} MDP has been proposed in \cite{Feinberg-Piunovskiy-2019}.  Loosely, this property means that the series of the queues of the entrance time to $\Delta$ go to zero uniformly over the set of policies. Recently, there has been  very active research on this topic trying to explore the links between, on one hand, the MDP being uniformly absorbing and, on the other hand, compactness of the set of occupation measures for a suitable topology; see, e.g., \cite{DuPri24-ESAIM,piunovskiy2023continuity,schal75}.

Regarding the topology on the set of occupation measures, there exist two main choices: the usual weak ($w$) topology and the weak-strong ($ws$) topology. They are both initial topologies on the set of finite measures on $\XX\times\AA$; they are the coarsest topologies making continuous the mappings 
$$\mu\mapsto \int_{\XX\times\AA} f(x,a)\mu(\mathrm{d}x,\mathrm{d}a),$$ where: for the weak topology we consider continuous and bounded functions $f:\XX\times\AA\rightarrow\RR$, while for the $ws$-topology we consider bounded Carath\'eodory functions (i.e., measurable functions $f:\XX\times\AA\rightarrow\RR$ which are continuous in $a\in\AA$ for every fixed $x\in\XX$). At this point, note that the $ws$-topology is finer than the weak topology. Also, concerning the continuity properties of the controlled Markov process, two main sets of hypotheses are usually considered in the literature: the strong (S) continuity conditions and the weak (W) continuity conditions which correspond, roughly, to the cases of a strong-Feller and a weak-Feller transition kernel, respectively.

The paper \cite{piunovskiy2023continuity} has studied compactness of the set of occupation measures for the weak topology under the (S) and the (W) conditions. It is proved that, for an isolated absorbing set $\Delta$ (i.e., $\Delta$ is both an open and a closed set), the MDP is uniformly absorbing if and only if the natural projection mapping, which associates to each strategic probability measure its corresponding occupation measure, is continuous. (Recall that the set of strategic probability measures is known to be a compact metric space \cite{balder89,schal75}.) In particular, the set of occupation measures of an MDP which is uniformly absorbing to an isolated set $\Delta$ is $w$-compact. The converse implication that an MDP with $w$-compact occupation measures is  uniformly absorbing is, however, mentioned in \cite{piunovskiy2023continuity} as an interesting open problem. 

The other reference dealing with this topic is \cite{DuPri24-ESAIM}. In that paper,  the authors deal with an absorbing MDP satisfying the (S) conditions and they consider the $ws$-topology on the set of occupation measures. The main result is that the MDP is uniformly absorbing if and only if the set of occupation measures is $ws$-compact. This result has been obtained without imposing any topological condition on the set $\Delta$ (in fact, the results in \cite{DuPri24-ESAIM} are proved for a measurable, not necessarily metric nor topological, state space $\XX$).

Even if the result in \cite{DuPri24-ESAIM} gives a complete overview of the topic when considering the (S) conditions and the $ws$-topology,
 there remain, however, several interesting open problems when considering the weak topology on the set of occupation measures. In addition to the implication mentioned in \cite{piunovskiy2023continuity}  whether an MDP with $w$-compact occupation measures is necessarily uniformly absorbing, it turns out that the condition that the absorbing set $\Delta$ is isolated is, in fact, fairly restrictive. Although, intuitively, one may think that assuming that the controlled process is absorbed by an isolated cemetery state can be done without loss of generality, this is far from being true. 
 
 The present paper is devoted, precisely, to explore the relations existing between three main issues arising in an absorbing MDP, namely, (i) $w$-compactness of the set of occupation measures; (ii) continuity of the natural projection mapping from the set of strategic probability measures onto the set of occupation measures; (iii) the MDP being uniformly absorbing. We study the various implications existing between these three main properties under both the (S) and the~(W) continuity conditions, and also based on the behavior of the dynamic system at the boundary of the absorbing set $\Delta$ or based on the topological properties of $\Delta$. 
 
 In this paper we propose, among others, conditions under which all statements (i)--(iii) are equivalent. The conditions that we propose here include, as a corollary, the case of an isolated absorbing set $\Delta$. Hence, not only we give a positive answer to the conjecture made in \cite{piunovskiy2023continuity}, but we also provide weaker conditions than $\Delta$ being isolated under which the full equivalence of (i)--(iii) holds. Interestingly enough, for an MDP satisfying the (S) conditions, we reach equivalence of (i)--(iii) without any topological condition on $\Delta$. Let us also point out that, while in \cite{piunovskiy2023continuity} the authors consider that all actions in $\AA$ are available at each state $x\in\XX$, here we allow for state-dependent action sets $\AA(x)$. This introduces an interesting gain of generality and, obviously, continuity and measurability issues on the multifunction $x\mapsto\AA(x)$ become necessary.

In the remaining of this section, we introduce some notation that we shall use throughout the paper. In Section \ref{sec-2} we give our main definitions and state some basic results.
Section \ref{sec-3} presents our main results for an absorbing MDP under the weak continuity conditions, while absorbing MDPs satisfying the strong optimality conditions are studied in  Section \ref{sec-4}.

\paragraph{Notation and terminology.}
The standard symbols $\NN=\{0,1,2,\ldots\}$ and  $\RR$ denote
the sets of natural and real numbers, respectively.
Given a measurable space $(\mathbf{\Omega},\bcal{F})$, lety $\bcal{M}^+(\mathbf{\Omega})$ and  $\bcal{P}(\mathbf{\Omega})$ be the sets of finite nonnegative measures and  probability measures, respectively, on $(\mathbf{\Omega},\bcal{F})$.
For a measurable set $\Gamma\in\bcal{F}$, we denote by $\mathbf{I}_{\Gamma}:\mathbf{\Omega}\rightarrow\{0,1\}$ the indicator function of the set~$\Gamma$.
Given $\omega\in\mathbf{\Omega}$, we write $\delta_{\omega}\in\bcal{P}(\mathbf{\Omega})$ for the Dirac probability measure concentrated at $\omega$. 

If $(\mathbf{\Omega}',{\bcal{F}'})$ is another measurable space, 
a (substochastic) kernel on $\mathbf{\Omega}'$ given $\mathbf{\Omega}$ is a mapping
$Q:\mathbf{\Omega}\times{\bcal{F}'}\rightarrow[0,1]$ such that $\omega\mapsto Q(B|\omega)$ is 
measurable on $(\mathbf{\Omega},\bcal{F})$ for every $B\in\bcal{F}'$,  and  such that $B\mapsto Q(B|\omega)$ is in 
$\bcal{M}^+(\mathbf{\Omega}')$
for every $\omega\in\mathbf{\Omega}$.  When $Q(\mathbf{\Omega}'|\omega)=1$ for every $\omega\in\mathbf{\Omega}$ we will  say that $Q$ is a \textit{stochastic kernel}.
Given $\Gamma\in\bcal{F}$, we denote by $\mathbb{I}_{\Gamma}$ the kernel on $\mathbf{\Omega}$ given $\mathbf{\Omega}$ defined by
$\mathbb{I}_{\Gamma}(B|\omega)=\mathbf{I}_{\Gamma}(\omega) \delta_{\omega}(B)$ for $\omega\in\mathbf{\Omega}$ and $B\in\bcal{F}$. In words, if $\omega\in\Gamma$ then $\mathbb{I}_\Gamma(dy|\omega)$ is the Dirac measure at $\omega$, and if $\omega\notin\Gamma$ then $\mathbb{I}_\Gamma(dy|\omega)$ is the null measure.

For $\mu\in\bcal{M}^{+}(\mathbf{\Omega})$ and a kernel $Q$  on  $\mathbf{\Omega}'$ given $\mathbf{\Omega}$, we define $\mu Q\in\bcal{M}^+(\mathbf{\Omega}')$ by means of  
$$B\mapsto \mu Q\,(B)= \int_{\mathbf{\Omega}} Q(B|\omega) \mu(\mathrm{d}\omega)\quad\hbox{for $B\in\bcal{F}'$}.$$
If $f:\mathbf{\Omega}'\rightarrow\RR$ is a measurable function, by $Qf$ we will denote the measurable function on $\mathbf{\Omega}$
$$Qf(\omega)= \int_{\mathbf{\Omega}'} f(\omega')Q(\mathrm{d}\omega'|\omega)\quad\hbox{for $\omega\in\mathbf{\Omega}$},$$
provided that the integrals are well defined.
If $({\mathbf{\Omega}'',}\bcal{F}'')$ is another measurable space and $R$ is a kernel on $\mathbf{\Omega}''$
given $\mathbf{\Omega}'$, then  $QR$ is the kernel on $\mathbf{\Omega}''$ given $\mathbf{\Omega}$ defined as
$$QR(\Gamma |\omega)= \int_{\mathbf{\Omega}'} R(\Gamma | \omega') Q(\mathrm{d}\omega' | \omega)
\quad\hbox{for $\Gamma\in\bcal{F}''$ and $\omega\in\mathbf{\Omega}$}.$$

Given the measurable product space $(\mathbf{\Omega}\times\mathbf{\Omega}',\bcal{F}\otimes\bcal{F}')$ and a measure $\mu\in\bcal{M}^+(\mathbf{\Omega}\times\mathbf{\Omega}')$, the marginal measures are denoted by 
$$\mu^{\mathbf{\Omega}}(\cdot)=\mu(\cdot\times \mathbf{\Omega}')\in\bcal{M}^+(\mathbf{\Omega})\quad\hbox{and}\quad
\mu^{\mathbf{\Omega}'}(\cdot)=\mu( \mathbf{\Omega}\times\cdot)\in\bcal{M}^+(\mathbf{\Omega}').
$$

On a topological space $\mathbf{S}$ we will always consider the Borel $\sigma$-algebra $\bcal{B}(\mathbf{S})$. 
The set of real-valued continuous functions defined on $\mathbf{S}$ will be denoted by $\mathbf{C}(\mathbf{S})$.
We will write $\mathbf{C}_{b}({\mathbf{S}})$ for the family of continuous and bounded continuous functions, and $\mathbf{C}_{b}^{+}({\mathbf{S}})$ for the set 
of non-negative bounded continuous functions on $\mathbf{S}$.

The space
$\bcal{M}^+(\mathbf{S})$ of finite measures on $(\mathbf{S},\bcal{B}(\mathbf{S}))$  will be endowed with
the weak  topology (which will be also referred to as the $w$-topology). It is the coarsest topology for which the mappings $\mu\mapsto \int f\mathrm{d}\mu$, for  
$f\in\mathbf{C}_{b}({\mathbf{S}})$, are continuous.
This topology is metrizable when $\mathbf{S}$ is separable. 
Convergence in the $w$-topology will be denoted by  $\mu_n\stackrel{w}{\rightarrow}\mu$.

\section{Definitions and preliminary results}
\label{sec-2}

\paragraph{The control model.}
The Markov control model is  $\bfrak{M}=(\mathbf{X},\mathbf{A},\{\mathbf{A}(x): x\in \mathbf{X}\},Q,\eta)$,
which consists of the following elements.
\begin{itemize}
\item The state space $\XX$ and the action space $\AA$, which are both Borel spaces endowed with their respective Borel  $\sigma$-algebras  $\bcal{B}(\XX)$ and  $\bcal{B}(\AA)$.
\item A family of nonempty measurable sets $\AA(x)\subseteq \AA$ for each $x\in\XX$, which stand for the set of actions available at $x$. We suppose that the set 
\[ \KK=\{(x,a)\in\XX\times\AA: a\in \AA(x)\}\]  of feasible state-action pairs belongs to $\bcal{B}(\XX)\otimes\bcal{B}(\AA)$ and that 
it contains the graph of a measurable function from $\XX$ to $\AA$.
\item A stochastic kernel $Q$ on $\mathbf{X}$ given $\KK$, which stands for the transition probability function. 
\item The initial distribution of the system $\eta\in\bcal{P}(\XX)$.
\end{itemize}

\paragraph{Control policies.}
The canonical space of the state-action process is defined as 
$\mathbf{\Omega}=(\XX\times\AA)^{\NN}$ and it is endowed with the product $\sigma$-algebra $\bcal{F}$. The coordinates projection functions are defined as follows: for $\omega\in\mathbf{\Omega}$ and $n\in\NN$,
\[ X_n(\omega)=x_n,\ A_n(\omega)=a_n,\ H_n(\omega)=(x_0,a_0,\ldots,x_n).\]
We will use the notation $\mathbf{H}_0=\XX$ and $\mathbf{H}_n=(\XX\times\AA)^{n}\times\XX$  for $n\ge1$. 

A control policy $\pi$  is a sequence $\{\pi_n\}_{n\in\NN}$ of stochastic kernels on $\AA$ given $\HH_n$ such that
$$\pi_n(\AA(x_n)|h_n)=1\quad\hbox{for every $n\in\NN$ and $h_n=(x_0,a_0,\ldots,x_n)\in\HH_n$}.$$  
The set of all policies is denoted by $\mathbf{\Pi}$.
The  set of stochastic kernels $\sigma$ on $\AA$ given $\XX$ satisfying $\sigma(\AA(x)|x)=1$
for any $x\in \XX$ is denoted  by $\mathbf{M}$.
Any $\sigma\in\mathbf{M}$ will be identified with the policy $\pi\in\mathbf{\Pi}$ given by
$$\pi_{n}(\cdot |h_{n})=\sigma(\cdot |x_{n})\quad\hbox{ for any $n\in\NN$ and $h_{n}=(x_{0},a_{0},\ldots,x_{n})\in \mathbf{H}_{n}$}.$$
In what follows, $\mathbf{M}$ will be referred to as the set of \emph{Markovian stationary policies}. These sets are nonempty as a consequence of the hypothesis that $\KK$ contains the graph of a measurable function.

\paragraph{Strategic probability measures.}
Given any policy $\pi \in \mathbf{\Pi}$, by the Ionescu-Tulcea theorem, there exists  a unique probability
measure~$\mathbb{P}_{\pi}\in \bcal{P}(\mathbf{\Omega})$
such that 
$\mathbb{P}_{\pi}(X_{0}\in B)=\eta(B)$, $$
\mathbb{P}_{\pi}(X_{n+1}\in B\mid H_{n},A_{n}) = Q(B\mid X_{n},A_{n}) \quad\hbox{and}\quad
\mathbb{P}_{\pi}(A_{n}\in C\mid H_{n})= \pi_{n}(C\mid H_{n})$$
 for every $n\in\NN$, $B\in \bcal{B}(\XX)$, and  $ C\in \bcal{B}(\mathbf{A})$,  with
$\mathbb{P}_{\pi}$-probability one.  We note that, by construction, we have that 
$$\mathbb{P}_\pi(H_n\in \KK^{n}\times\XX)=1\quad\hbox{for each $n\in\NN$},$$  with the convention that $\KK^0\times\XX=\XX$, and also that 
$ \mathbb{P}_{\pi} (\KK^{\NN})=1$. 
The expectation operator associated to  $\mathbb{P}_{\pi}$ will be denoted by $\mathbb{E}_{\pi}$ (note that this is a slight abuse of notation since, strictly speaking, we should write $\mathbb{E}_{\mathbb{P}_\pi}$).
We will refer to $\mathbb{P}_\pi$ as to a \emph{strategic probability measure}, and we will denote by
 \[\bcal{S}=\{\mathbb{P}_\pi\}_{\pi\in\mathbf{\Pi}}\subseteq\bcal{P}(\mathbf{\Omega})\]
the set of all strategic probability measures. 

\paragraph{Topologies on $\bcal{P}(\mathbf{\Omega})$.}
On the set  $\bcal{P}(\mathbf{\Omega})$ of probability measures on $(\mathbf{\Omega},\bcal{F})$ we will consider two topologies. 
The $ws^\infty$-topology is the coarsest topology making continuous
the mappings
\[ \mathbb{P}\mapsto \int_{\mathbf{H}_n} f(x_0,a_0,\ldots,x_n)\mathrm{d}(\mathbb{P}\circ H_n^{-1})\]
for every $n\in\NN$ and  every bounded Carath\'eodory  function $f:\mathbf{H}_n\rightarrow\RR$, i.e., $f$ is measurable and $(a_0,\ldots,a_{n-1})\mapsto f(x_0,a_0,\ldots,a_{n-1},x_n)$ is continuous on $\AA^n$ for every fixed $(x_0,\ldots,x_n)\in\XX^{n+1}$; see section 4 in \cite{schal75} for more details on this topology.
The $ws^\infty$-convergence on~$\bcal{P}(\mathbf{\Omega})$ will be denoted by $\mathbb{P}_n\stackrel{ws^{\infty}}{\longrightarrow}\mathbb{P}$.

We can also consider the usual weak topology (or $w$-topology), which  is metrizable (indeed, $\mathbf{\Omega}$ is the countable product of Borel spaces and it is itself a Borel space with the product topology). The $w$-topology can be also defined as the coarsest topology making continuous the  mappings
\begin{equation}\label{eq-carath}
 \mathbb{P}\mapsto \int_{\mathbf{H}_n} f(x_0,a_0,\ldots,x_n)\mathrm{d}(\mathbb{P}\circ H_n^{-1})
 \end{equation}
for every $n\in\NN$ and  every $f\in\mathbf{C}_b(\mathbf{H}_n)$; see Lemma 4.1 in \cite{schal75}.
Convergence on~$\bcal{P}(\mathbf{\Omega})$ for this topology will be denoted by $\mathbb{P}_n\stackrel{w}{\rightarrow}\mathbb{P}$.
It is clear that $\mathbb{P}_n\stackrel{ws^{\infty}}{\longrightarrow}\mathbb{P}$ implies $\mathbb{P}_n\stackrel{w}{\rightarrow}\mathbb{P}$ and so the $w$-topology is coarser than the $ws^\infty$-topology.  

\paragraph{Absorbing control models.} Let $\Delta\in\bcal{B}(\XX)$ be some nonempty measurable subset of the state space. 
The hitting time $T_\Delta$ of the set $\Delta$ is given by  $T_\Delta:\mathbf{\Omega}\rightarrow\NN\cup\{\infty\}$ defined as
$$T_\Delta(x_0,a_0,x_1,a_1,\ldots)=\min\{n\ge0: x_n\in\Delta\},$$
where the $\min$ over the empty set is defined as $+\infty$. 
Next we give the definition of  absorbing and uniformly absorbing control models.

\begin{definition}\label{def-absorbing}
\label{absorbing} Given the control model $\bfrak{M}$ and  $\Delta\in\bcal{B}(\XX)$, 
we say that 
$\bfrak{M}$ is absorbing  to $\Delta$  if the conditions (a)--(b) below are satisfied, and we say that 
$\bfrak{M}$ is uniformly absorbing  to $\Delta$ if, additionally, condition (c) holds.
\begin{enumerate}[label=(\alph*).]
\item 
For every $(x,a)\in\KK$ such that $x\in\Delta$ we have $Q(\Delta|x,a)=1$.
\item For any $\pi\in\mathbf{\Pi}$ the expected hitting time $\mathbb{E}_{\pi}[T_\Delta]$ is finite.
\item The following limit holds: $$\lim_{n\rightarrow\infty} \sup_{\pi\in\mathbf{\Pi}} \sum_{t>n}\mathbb{P}_{\pi}\{T_\Delta>t\}=0.$$
\end{enumerate}
\end{definition}
As noted in \cite{Feinberg-Piunovskiy-2019} ---see the comments after Definition 3.6 therein--- the supremum in (c) can be equivalently taken over the set $\mathbf{M}$ of Markovian stationary policies.

\begin{definition}\label{def-occupation-measure}
If the control model $\bfrak{M}$ is absorbing to $\Delta$, then we define the occupation measure $\mu_\pi\in\bcal{M}^+(\XX\times\AA)$ 
 of the control policy $\pi\in\mathbf{\Pi}$ as
$$
\mu_{\pi}(\Gamma) = \mathbb{E}_{\pi}\Big[ \sum_{t=0}^\infty \mathbf{I}_{\{T_\Delta>t\}}\cdot\mathbf{I}_{\Gamma}(X_t,A_{t})\Big] 
= \sum_{t=0}^\infty \mathbb{P}_{\pi}\{X_t\in\Delta^c, (X_t,A_t)\in\Gamma\} \quad\hbox{for $\Gamma\in\bcal{B}(\XX)\otimes\bcal{B}(\AA)$.}
$$
The set of all occupation measures is
$\bcal{O}=\{\mu_{\pi}:\pi\in\mathbf{\Pi}\}\subseteq\MM^+(\XX\times\AA)$.
\end{definition}

Roughly, the occupation measure counts the expected number of visits of the state-action process   to the set  $\Gamma$ up to time~$T_\Delta$   and it vanishes on $\Delta$, that is, $\mu_{\pi}(\Delta\times\AA)=0$. We also have that $\mu_\pi(\KK^c)=0$. 
 It is easy to see that 
\[\mu_\pi(\XX\times\AA)=\sum_{t=0}^\infty\mathbb{P}_{\pi}\{T_\Delta>t\}=\mathbb{E}_\pi[T_\Delta]<\infty\] 
(recall Definition \ref{def-absorbing}(b)) and so $\mu_\pi\in\MM^+(\XX\times\AA)$ for every $\pi\in\mathbf{\Pi}$. As a consequence of \cite[Sections~4.4 and~5.5]{dynkin79},  Definition \ref{def-absorbing}(b) implies that, in fact,  $\sup_{\pi\in\mathbf{\Pi}}\mathbb{E}_\pi[T_\Delta]$ is finite. 
Therefore, the set of occupation measures
is bounded  in  $ \MM^+(\XX\times\AA)$, meaning that $\sup_{\mu\in\bcal{O}} \mu(\XX\times\AA)<\infty$.

\paragraph{The projection mapping.}
Let $\bfrak{M}$ be absorbing to $\Delta$. 
There is a natural projection
 $\bfrak{P}:\bcal{S}\rightarrow\bcal{O}$ from the set of strategic probability measures to the set of occupation measures defined as follows. If $\mathbb{P}\in\bcal{S}$ then $\mu=\bfrak{P}(\mathbb{P})$ is given by 
 $$\mu(\Gamma)=
 \sum_{t=0}^\infty \mathbb{P}\{X_t\in\Delta^c, (X_t,A_t)\in\Gamma\} \quad\hbox{for $\Gamma\in\bcal{B}(\XX)\otimes\bcal{B}(\AA)$.}$$
  It is clear that when $\mathbb{P}$ is the strategic probability measure of a policy $\pi\in\mathbf{\Pi}$, that is, $\mathbb{P}=\mathbb{P}_\pi$, then we have
 $\bfrak{P}(\mathbb{P}_\pi)=\mu_\pi$.

\paragraph{The characteristic equations.} In what follows, assume that $\bfrak{M}$ is absorbing to $\Delta$. 
We say that a measure $\mu\in\MM^+(\XX\times\AA)$ satisfies the \emph{characteristic equations} when 
\begin{equation}\label{eq-linear-equation} 
\mu(\KK^c)=0\quad\hbox{and}\quad \mu^\XX=(\eta+\mu Q)\mathbb{I}_{\Delta^c}.
\end{equation}
Note that the last equation can be written, in a less concise way, as 
\begin{equation}\label{eq-char-eq}
\mu^\XX(B)=\eta(B\cap\Delta^c)+\int_{\XX\times\AA} Q(B\cap \Delta^c|x,a)\mu(\mathrm{d}x,\mathrm{d}a)\quad\hbox{for $B\in\bcal{B}(\XX)$}.
\end{equation}
Since a measure $\mu$ satisfying the characteristic equations is supported on $\KK$, it is not relevant that $Q(B\cap\Delta^c|x,a)$ in \eqref{eq-char-eq} is not defined when $(x,a)$ is outside $\KK$.
We denote by $\bcal{C}$ the set of measures in $\MM^+(\XX\times\AA)$ which satisfy the characteristic equations. 

We state the following well-known fact. For a proof, we refer to Proposition 3.3(i) in \cite{DuPri24-ESAIM}.

\begin{lemma} Let $\bfrak{M}$ be absorbing to $\Delta$. Any occupation measure satisfies the characteristic equations, that is, 
$\bcal{O}\subseteq\bcal{C}$.
\end{lemma}

Our Lemma \ref{lemma-1} below requires a disintegration result of a finite measure supported on $\KK$ by means of a Markov policy in $\mathbf{M}$. To this end, we need to impose  that the action sets $\AA(x)$ are closed  and that the correspondence $x\mapsto\AA(x)$ is weakly measurable. In particular,  by \cite[Theorems 18.6 and 18.13]{aliprantis06}, this will imply that $\KK$ is measurable and that it
it contains the graph of a measurable function from $\XX$ to $\AA$ (cf. the definition of the control model).

\begin{lemma}
\label{lemma-disintegration} 
Suppose that the action sets $\AA(x)$ are closed for every $x\in\XX$ and that the correspondence $x\mapsto\AA(x)$ is weakly measurable.
  For every  $\mu\in\bcal{M}^{+}(\XX\times\AA)$ 
such that $\mu(\mathbf{K}^{c})=0$
there exists $\sigma\in\mathbf{M}$ 
with
$\mu= \mu^{\XX}\otimes \sigma$.

There is uniqueness in the following sense: if $\sigma,\sigma'\in\mathbf{M}$  satisfy $\mu= \mu^{\XX}\otimes \sigma= \mu^{\XX}\otimes \sigma'$ then, $\mu^\XX$-almost surely,
we have $\sigma(\mathrm{d}a|x)=\sigma'(\mathrm{d}a|x)$.
\end{lemma}

Our next result is borrowed from Proposition 3.8(i) in \cite{DuPri24-ESAIM} whose proof relies, precisely, on the disintegration Lemma \ref{lemma-disintegration}.

\begin{lemma}\label{lemma-1}
Let the control model $\bfrak{M}$ be absorbing to $\Delta$, and suppose that the action sets $\AA(x)$ are closed for every  $x\in\XX$ and that the correspondence $x\mapsto\AA(x)$ is weakly measurable. If two occupation measures
If $\mu_{\pi},\mu_{\pi'}\in\bcal{O}$ are such that $\mu_{\pi}-\mu_{\pi'}\in\bcal{M}^+(\XX\times\AA)$ then $\mu_{\pi}=\mu_{\pi'}$. 
\end{lemma}
\textit{Proof.} 
We have that  $\nu=\mu_{\pi}-\mu_{\pi'}$ is an invariant  measure for the kernel $Q\mathbb{I}_{\Delta^c}$ as introduced in  \cite[Definition 3.6]{DuPri24-ESAIM}.
Proceeding by contradiction, if $\nu$ is not the null measure then, using \cite[Proposition 3.8(i)]{DuPri24-ESAIM}, we derive that $\nu^\XX\perp\mu_\pi^\XX$. However, from the equation $\nu+\mu_{\pi'}=\mu_{\pi}$ we obtain $\nu^\XX\ll\mu_\pi^\XX$. This is a contradiction and we conclude that  $\nu=\mu_{\pi}-\mu_{\pi'}$ is the null measure.
\hfill$\Box$
\\[10pt]
\indent
This yields the  interesting fact that, given two occupation measures $\mu_{\pi},\mu_{\pi'}\in\bcal{O}$, it is not possible that one of them dominates the other one (i.e., either $\mu_{\pi}-\mu_{\pi'}$ or $\mu_{\pi'}-\mu_\pi$ are nonnegative measures) unless they are equal.

\paragraph{Continuity conditions.}
In the sequel we shall impose two different sets of conditions on our control model. (Note that we incorporate the weak measurability of the actions sets correspondence to these conditions so that we will be in position to use Lemmas \ref{lemma-disintegration} and \ref{lemma-1}.) The so-called \emph{weak continuity conditions} are the following.
\\[10pt]
\noindent\textbf{Conditions (W)}
\begin{conditW}
\item \label{HypW-multifunct}  The multifunction $x\mapsto\AA(x)$ is compact-valued, upper semi-continuous, and weakly measurable.
\item \label{HypW-transition} The mapping $(x,a)\mapsto \int_\XX f(y)Q(\mathrm{d}y|x,a)=Qf(x,a)$ is continuous on $\KK$ for any $f\in\mathbf{C}_b(\XX)$. 
\end{conditW}

Next we state the \emph{strong continuity conditions}
\\[10pt]
\noindent\textbf{Conditions (S)}
\begin{conditS}
\item \label{S1} The action set $\AA(x)$ is compact for every $x\in\XX$  and  $x\mapsto\AA(x)$ is  weakly measurable.
\item \label{S2} The mapping $a\mapsto Q(B|x,a)$ is continuous on $\AA(x)$ for each fixed $B\in\bcal{B}(\XX)$ and $x\in\XX$.
\end{conditS}

These continuity conditions are fairly standard and they have been widely used in the literature; see, for instance,  \cite{schal75}. The words weak and strong should not be misleading since none the above conditions  imply each other: the conditions that they impose are, in fact, of different nature.
A very important consequence of these continuity conditions is the following.
 For a proof of this result we refer to \cite{balder89}.

\begin{lemma}\label{lem-balder} 
If the control model $\bfrak{M}$ satisfies either Conditions (W) or (S), then $\bcal{S}$ is a compact metric space for the $w$-topology. Besides, under the Conditions (S), the $ws^\infty$-topology and the $w$-topology on $\bcal{S}$ coincide.
\end{lemma}

Throughout this paper we will discuss the relation existing between the following \emph{three basic properties} of an absorbing control model $\bfrak{M}$:
\begin{itemize}
\item The set $\bcal{O}$ of occupation measures is $w$-compact,
\item The control model $\bfrak{M}$ is uniformly absorbing to $\Delta$.
\item The natural projection mapping $\bfrak{P}:\bcal{S}\rightarrow\bcal{O}$ is continuous (here, $\bcal{S}$ and $\bcal{O}$ are endowed with their respective weak topologies).
\end{itemize}

In view of Lemma \ref{lem-balder} and under the Conditions (W) or (S), there is an implication which is always true: if $\bfrak{P}$ is continuous then $\bcal{O}$ is $w$-compact. 
To obtain further implications between these properties, some additional conditions will be needed. 
But, first of all, let us  give a result which holds, in general, under the Conditions (W) or (S).

\begin{proposition}\label{prop-ua-rel-compact} 
Suppose that $\bfrak{M}$ satisfies the Conditions (W) or (S) and that it is uniformly absorbing to $\Delta$. Then $\bcal{O}$ is relatively $w$-compact.
\end{proposition}
\textit{Proof.}
We will use the following characterization of relative $w$-compactness taken from \cite[Theorem 4.5.10]{bogachev18}. The set $\bcal{O}$, which is bounded in $\bcal{M}^+(\XX\times\AA)$, is relatively $w$-compact if and only if for any sequence $\{f_n\}$ in $\mathbf{C}^+_b(\XX\times\AA)$  such that $f_n\downarrow0$  we have 
$$\lim_{n\rightarrow\infty} \sup_{\mu\in\bcal{O}}\int_{\XX\times\AA} f_n\mathrm{d}\mu=0.$$

Suppose now that $\mathfrak{M}$ is uniformly absorbing to $\Delta$ and let $\{f_n\}$ be a sequence of functions satisfying the above conditions. Let the constant $\mathbf{c}>0$ be such that $0\le f_n\le \mathbf{c}$ for every~$n$. Given $\epsilon>0$ there exists some $k_0$ such that 
$$ \sup_{\pi\in\mathbf{\Pi}} \sum_{t>k_0}\mathbb{P}_{\pi}\{T_\Delta>t\}\le\epsilon/2\mathbf{c}.$$
Given any $\mu\in\bcal{O}$, let   $\pi\in\mathbf{\Pi}$ be such that $\mu=\mu_\pi$. We have 
\begin{eqnarray*}
\int_{\XX\times\AA} f_n(x,a)\mu(\mathrm{d}x,\mathrm{d}a)&=&\mathbb{E}_\pi\Big[\sum_{t=0}^\infty f_n(X_t,A_t)\mathbf{I}_{\{T_\Delta>t\}} \Big]\\
&\le&  \mathbb{E}_\pi\Big[\sum_{t=0}^{k_0} f_n(X_t,A_t)\Big]+\mathbf{c}\sum_{t>{k_0}}\mathbb{P}_\pi\{T_\Delta>t\}\\
&\le& \mathbb{E}_\pi\Big[\sum_{t=0}^{k_0} f_n(X_t,A_t)\Big]+\epsilon/2.
\end{eqnarray*}
Consequently,
\begin{eqnarray}
\sup_{\mu\in\bcal{O}} \int_{\XX\times\AA} f_n\mathrm{d}\mu 
&\le& \sup_{\pi\in\mathbf{\Pi}} \mathbb{E}_\pi\Big[\sum_{t=0}^{k_0} f_n(X_t,A_t)\Big]+\epsilon/2 \label{eq-rel-w-compact} \\
&=& \sup_{\pi\in\mathbf{\Pi}}\int_{\mathbf{\Omega}}\Big[ \sum_{t=0}^{k_0} f_n(x_t,a_t)\Big]\mathrm{d}\mathbb{P}_\pi  +\epsilon/2. \nonumber
\end{eqnarray}
We have that  $\sum_{t=0}^{k_0} f_n(x_t,a_t)$ is a sequence (in $n$) of bounded and continuous functions which decrease to $0$. Dominated or monotone convergence yields
\begin{equation}\label{eq-dini}
\int_{\mathbf{\Omega}}\Big[\sum_{t=0}^{k_0} f_n(x_t,a_t)\Big]\mathrm{d}\mathbb{P}_\pi\downarrow0\quad\hbox{ as $n\rightarrow\infty$ for every $\pi\in\mathbf{\Pi}$}.
\end{equation}
However, by definition of the $w$-topology, we have that the mapping
$$\mathbb{P}\mapsto \int_{\mathbf{\Omega}}\Big[\sum_{t=0}^{k_0} f_n(x_t,a_t)\Big]\mathrm{d}\mathbb{P}\quad\hbox{ is continuous on $\bcal{P}(\mathbf{\Omega})$}.$$
It follows from Dini's theorem that the convergence in \eqref{eq-dini} is uniform on the compact set $\bcal{S}$ (recall Lemma \ref{lem-balder}).
Hence, there exists some $n_0$ such that $n\ge n_0$ implies 
$$\sup_{\pi\in\mathbf{\Pi}} \mathbb{E}_\pi\Big[\sum_{t=0}^{k_0} f_n(X_t,A_t)\Big]\le \epsilon/2.$$
The result now follows from \eqref{eq-rel-w-compact}. This concludes the proof that $\bcal{O}$ is relatively $w$-compact when~$\bfrak{M}$ is uniformly absorbing.\hfill$\Box$\\[10pt]
\indent
In the following sections, we will put ourselves in the context of either the weak or the strong continuity conditions.

\section{Weakly continuous absorbing control models}\label{sec-3}

Throughout this section we will assume that the control model $\bfrak{M}$ is absorbing to $\Delta$ and that it satisfies the Conditions (W). 
As a direct consequence of the Condition \ref{HypW-multifunct} and using \cite[Theorem 17.10]{aliprantis06}, we derive  that $\KK$ is a closed subset of $\XX\times\AA$ since it is the graph of a compact-valued upper semicontinuous multifunction on the metric space $\AA$.

Next we develop several conditions yielding some implications between $w$-compactness of $\bcal{O}$, continuity of $\bfrak{P}$, and $\bfrak{M}$ being uniformly absorbing.

\subsection{Strongly continuous boundary condition at $\Delta$}

We introduce the next condition.

\begin{condition}[\textbf{Condition (SCB)}]
Given a control model $\bfrak{M}$ absorbing to $\Delta$, we will say that it satisfies the strongly continuous boundary (SCB) condition at $\Delta$ when, for every measurable $B\subseteq\partial\Delta$, the function $(x,a)\mapsto Q(B|x,a)$ is continuous on $\KK$.
\end{condition}

In the sequel, we shall just mention that $\bfrak{M}$ satisfies the Condition (SCB) and we will omit that it is satisfied at $\Delta$. 
Our next result gives some important consequences of this condition. In particular, item (ii) yields a weak continuity-like result for the transition kernel when the term $\mathbf{I}_{\Delta^c}$ is incorporated (cf. Condition \ref{HypW-transition}).

\begin{lemma}\label{lem-cont-boundary} 
Suppose that  $\bfrak{M}$ is absorbing to $\Delta$ and that it satisfies  the Conditions (W)   and (SCB).
\begin{description} 
\item[(i).] 
For any bounded measurable function $f:\XX\rightarrow\RR$, the mapping 
$$(x,a)\mapsto \int_{\partial{\Delta}} f(y)Q(\mathrm{d}y|x,a) \ \hbox{
is continuous on $\KK$.}$$
\item[(ii).] 
For any  $f\in\mathbf{C}_b(\XX)$, the mapping 
 $$(x,a)\mapsto\int_\XX f(y)\mathbf{I}_{\Delta^c}(y) Q(\mathrm{d}y|x,a)\ \hbox{
is continuous on $\KK$.}$$
\item[(iii).] The set $\bcal{C}$ is $w$-closed.
 \end{description}
\end{lemma}
\textit{Proof.} (i). Since the function $f$ is bounded, it is the increasing limit $f_n\uparrow f$ (resp. decreasing limit $f_n\downarrow f$) of simple bounded functions $\{f_n\}$. As a consequence of the Condition (SCB), for each $n\in\NN$, the function $(x,a)\mapsto \int_{\partial{\Delta}} f_n(y)Q(\mathrm{d}y|x,a)$ is continuous on $\KK$. We deduce that $(x,a)\mapsto \int_{\partial{\Delta}} f(y)Q(\mathrm{d}y|x,a)$ is a lower  (resp. upper) semicontinuous function on $\KK$. The proof of this statement is now complete.

(ii). Suppose first that $f$ is a nonnegative function. We can write for any $(x,a)\in\KK$
\begin{eqnarray}\label{eq-lemma-boundary}
\int_\XX f(y)\mathbf{I}_{\Delta^c}(y) Q(\mathrm{d}y|x,a) &=& 
\int_\XX f(y)\mathbf{I}_{\overline{\Delta^c}}(y) Q(\mathrm{d}y|x,a)-\int_\XX f(y)\mathbf{I}_{\Delta\cap\partial\Delta}(y) Q(\mathrm{d}y|x,a).
\end{eqnarray}
The function $y\mapsto f(y)\mathbf{I}_{\overline{\Delta^c}}(y)$ is the product of  function in $\mathbf{C}_b^+(\XX)$  and a bounded nonnegative upper semicontinuous function (since it is the indicator function of a closed set). Hence, $$(x,a)\mapsto \int_\XX f(y)\mathbf{I}_{\overline{\Delta^c}}(y) Q(\mathrm{d}y|x,a)\ \hbox{is upper semicontinuous on $\KK$.}$$ The rightmost term in \eqref{eq-lemma-boundary} is a continuous function of $(x,a)$ as a consequence of item (i). We conclude that the function in \eqref{eq-lemma-boundary} is upper semicontinuous. Observe now that we also have
\begin{eqnarray*}
\int_\XX f(y)\mathbf{I}_{\Delta^c}(y) Q(\mathrm{d}y|x,a) &=& 
\int_\XX f(y)\mathbf{I}_{(\Delta^c)^\circ}(y) Q(\mathrm{d}y|x,a)+\int_\XX f(y)\mathbf{I}_{\Delta^c\cap\partial\Delta}(y) Q(\mathrm{d}y|x,a).
\end{eqnarray*}
Reasoning similarly, we obtain that the lefthand side function is lower semicontinuous (this time, $(x,a)\mapsto\int_\XX f(y)\mathbf{I}_{(\Delta^c)^\circ}(y) Q(\mathrm{d}y|x,a)$ is lower semicontinuous). We conclude that the mapping $$(x,a)\mapsto\int_\XX f(y)\mathbf{I}_{\Delta^c}(y) Q(\mathrm{d}y|x,a)$$ is indeed  continuous  on $\KK$ whenever $f$ is nonnegative. The result for a general $f$ just follows by considering the positive and negative parts of $f$. 

(iii). 
To show that  $\bcal{C}$ is $w$-closed, assume that $\{\mu_n\}_{n\in\NN}$ is a sequence in $\bcal{C}$ that converges in the weak topology to some measure $\mu\in\bcal{M}^+(\XX\times\AA)$.
Choose an arbitrary 
$f\in\mathbf{C}_b^+(\XX)$ taking values in, e.g., $[0,\mathbf{c}]$. We know from item (ii) that  $Qf\mathbf{I}_{\Delta^c}$ is a continuous function on $\KK$ which takes values in $[0,\mathbf{c}]$. Since $\KK$ is closed, using Urysohn's lemma  \cite[Lemma 2.46]{aliprantis06}, we can extend $Qf\mathbf{I}_{\Delta^c}$ to a continuous function $h:\XX\times\AA\rightarrow[0,\mathbf{c}]$. 

Observe now that  since $\mu_n(\KK^c)=0$ for every $n\in\NN$ and $\KK$ is closed, we also have $\mu(\KK^c)=0$.  For every $n\in\NN$ 
\begin{eqnarray*}
\int_\XX f \mathrm{d} \mu_{n}^\XX &=&\int_\XX f\mathbf{I}_{\Delta^c} \mathrm{d}\eta+
\int_{\XX\times\AA}\int_\XX f(y)\mathbf{I}_{\Delta^c}(y) Q( \mathrm{d}y|x,a)\mu_{n}( \mathrm{d}x, \mathrm{d}a)\\
&=& \int_\XX f\mathbf{I}_{\Delta^c} \mathrm{d}\eta+
\int_{\XX\times\AA} h(x,a)\mu_{n}( \mathrm{d}x, \mathrm{d}a).
\end{eqnarray*}
Recalling that $h\in\mathbf{C}_b^+(\XX\times\AA)$ we can take the limit as $n\rightarrow\infty$ and obtain
\begin{eqnarray*}
\int_\XX f \mathrm{d} \mu^\XX &=&  \int_\XX f\mathbf{I}_{\Delta^c} \mathrm{d}\eta+
\int_{\XX\times\AA} h(x,a)\mu( \mathrm{d}x, \mathrm{d}a) \\
&=&
\int_\XX f\mathbf{I}_{\Delta^c} \mathrm{d}\eta+
\int_{\XX\times\AA}\int_\XX f(y)\mathbf{I}_{\Delta^c}(y) Q( \mathrm{d}y|x,a)\mu( \mathrm{d}x, \mathrm{d}a).\end{eqnarray*}
Since this equation holds for any $f\in\mathbf{C}_b^+(\XX)$, we conclude that $\mu$ indeed satisfies the characteristic equations, that is, $\mu\in\bcal{C}$.
\hfill$\Box$
\\[10pt]\indent
We state our main result under the Conditions (W) and (SCB). 

\begin{theorem}\label{th-W-SCB}
Suppose that the control model $\bfrak{M}$ is absorbing to $\Delta$ and that it satisfies the Conditions (W) and (SCB).
The following implications hold.
$$ \hbox{$\bfrak{P}$ is continuous} \ \Rightarrow\  \hbox{$\bfrak{M}$ is uniformly absorbing}\ \Rightarrow\ 
\hbox{$\bcal{O}$ is $w$-compact}.$$
\end{theorem}
\textit{Proof of the first implication.} We observe the following fact. Since, obviously, the function which associates to each $\mu_\pi\in\bcal{O}$ its total mass $$\mu_\pi(\XX\times\AA)=\int_{\XX\times\AA} 1\mu_\pi(\mathrm{d}x,\mathrm{d}a)=\mathbb{E}_{\pi}[T_\Delta]$$ is continuous for the weak topology on $\bcal{O}$, it turns out that when $\{\pi_n\}$ and $\pi^*$ in $\mathbf{\Pi}$ are such that $\mathbb{P}_{\pi_n}\stackrel{w}{\rightarrow}\mathbb{P}_{\pi^*}$ then we also have $\mathbb{E}_{\pi_n}[T_\Delta]\rightarrow \mathbb{E}_{\pi^*}[T_\Delta]$, just by composition of two continuous functions.

We establish, first of all, some useful identities. For every policy $\pi\in\mathbf{\Pi}$ and any $k\ge1$ 
\begin{eqnarray*}
\mathbb{E}_{\pi} [T_\Delta] &=& \sum_{t=0}^{k}\mathbb{P}_{\pi}\{X_t\in \Delta^c\} + \sum_{t>k}   \mathbb{P}_{\pi}\{T_\Delta>t\}\\
&=& \eta(\Delta^c)+\sum_{t=1}^{k} \mathbb{E}_{\pi}[Q(\Delta^c|X_{t-1},A_{t-1})] + \sum_{t>k}   \mathbb{P}_{\pi}\{T_\Delta>t\}.
\end{eqnarray*}
Using Urysohn's lemma \cite[Lemma 2.46]{aliprantis06} again, we can extend $(x,a)\mapsto Q(\Delta^c|x,a)$ which is continuous on $\KK$ (indeed, just let $f\equiv1$ in Lemma \ref{lem-cont-boundary}(ii)) to a continuous function $h:\XX\times\AA\rightarrow[0,1]$. We thus have
\begin{equation}\label{eq-tool-dugundji}
\mathbb{E}_{\pi} [T_\Delta] 
= \eta(\Delta^c)+\sum_{t=1}^{k} \mathbb{E}_{\pi}[h(X_{t-1},A_{t-1})] + \sum_{t>k}   \mathbb{P}_{\pi}\{T_\Delta>t\}
\quad\hbox{for any $\pi\in\mathbf{\Pi}$ and $k\ge1$}.
\end{equation}
The proof proceeds by contradiction. Hence, if $\bfrak{M}$ is not uniformly absorbing to $\Delta$ then there exist $\epsilon>0$ and sequences $k_n\uparrow\infty$ and $\{\pi_n\}$ in $\mathbf{\Pi}$ such that 
$$ \sum_{t>k_n} \mathbb{P}_{\pi_n}\{T_\Delta>t\}\ge\epsilon\quad\hbox{for every $n\in\NN$}.$$
Without loss of generality, we can assume the existence of some $\pi^*\in\mathbf{\Pi}$ such that $\mathbb{P}_{\pi_n}\stackrel{w}{\rightarrow}\mathbb{P}_{\pi^*}$ with,  as already mentioned, $\mathbb{E}_{\pi_n}[T_\Delta]\rightarrow \mathbb{E}_{\pi^*}[T_\Delta]$. For the policy $\pi^*$ we can find some $k^*$ such that
$$ \sum_{t>k^*} \mathbb{P}_{\pi^*}\{T_\Delta>t\}\le\epsilon/2.$$
For every $n$ such that $k_n\ge k^*$ we have $\sum_{t>k^*} \mathbb{P}_{\pi_n}\{T_\Delta>t\}\ge\epsilon$
and thus from  \eqref{eq-tool-dugundji} we get that
\begin{eqnarray*}
\epsilon/2 &\le & \sum_{t>k^*} \mathbb{P}_{\pi_n}\{T_\Delta>t\}-\sum_{t>k^*} \mathbb{P}_{\pi^*}\{T_\Delta>t\}\\
&=&\mathbb{E}_{\pi_n} [T_\Delta]-\mathbb{E}_{\pi^*} [T_\Delta]-
\sum_{t=1}^{k^*}\big( \mathbb{E}_{\pi_n}[h(X_{t-1},A_{t-1})]- \mathbb{E}_{\pi^*}[h(X_{t-1},A_{t-1})]\big).
\end{eqnarray*}
This is in contradiction with $\mathbb{E}_{\pi_n}[T_\Delta]\rightarrow \mathbb{E}_{\pi^*}[T_\Delta]$ and the fact that
\begin{equation}\label{eq-tool-finite-dim}
\mathbb{E}_{\pi_n}[h(X_{t-1},A_{t-1})] \rightarrow
 \mathbb{E}_{\pi^*}[h(X_{t-1},A_{t-1})]\quad\hbox{for any $t\ge1$}
\end{equation}
because $h\in\mathbf{C}_b(\XX\times\AA)$. We have thus established the first implication.
\\[5pt]\noindent
\textit{Proof of the second implication.}
 We know from Proposition \ref{prop-ua-rel-compact}  that $\bcal{O}$ is relatively $w$-compact. To obtain compactness, we must prove that $\bcal{O}$ is $w$-closed. 
 Thus,
consider a sequence $\{\pi_n\}_{n\in\NN}$ in $\mathbf{\Pi}$ such that $\mu_{\pi_n}\stackrel{w}{\rightarrow}\mu$. We know from Lemma \ref{lem-cont-boundary}(iii) that  
$\mu\in\bcal{C}$. We must prove that $\mu\in\bcal{O}$, that is, $\mu$ is the occupation measure of some policy.

By Lemma \ref{lem-balder}, we can  assume the existence of some $\pi^*\in\mathbf{\Pi}$ such that $\mathbb{P}_{\pi_n}\stackrel{w}{\rightarrow}\mathbb{P}_{\pi^*}$. Our first step in this proof is to show that 
$\mu^\XX=\mu_{\pi^*}^\XX$.
Let $f:\XX\rightarrow[0,1]$ be a continuous function. As we  have $\mu_{\pi_n}^\XX\stackrel{w}{\rightarrow}\mu^\XX$, for any $\epsilon>0$ there is some $n_0$ such that for every $n\ge n_0$
$$\int_\XX f\mathrm{d}\mu^\XX-\epsilon\le \int_\XX f\mathrm{d}\mu^\XX_{\pi_n}.$$ 
Since $\bfrak{M}$ is uniformly absorbing, for such $\epsilon>0$ there exists $k$ such that
$$
\sup_{\pi\in\mathbf{\Pi}}\sum_{t>k}^\infty \mathbb{P}_{\pi}\{T_\Delta>t\} \le\epsilon.
$$
For all $n\ge n_0$ we have 
\begin{equation}\label{eq-cond-UA-proof}
\int_\XX f\mathrm{d}\mu^\XX -\epsilon \le \int_\XX f\mathrm{d}\mu^\XX_{\pi_n}=\mathbb{E}_{\pi_n} \Big[ \sum_{t=0}^\infty  f(X_t)\mathbf{I}_{\Delta^c}(X_t)\Big]
\le  \mathbb{E}_{\pi_n} \Big[ \sum_{t=0}^k f(X_t)\mathbf{I}_{\Delta^c}(X_t)\Big]+\epsilon. 
\end{equation}
By Lemma \ref{lem-cont-boundary}(ii), the function 
 $$(x,a)\mapsto\int_\XX f(y)\mathbf{I}_{\Delta^c}(y) Q(\mathrm{d}y|x,a)$$ is continuous on $\KK$ and it can be extended to a continuous function $h$ on $\XX\times\AA$ taking values in $[0,1]$.
 Therefore, for $t\ge1$,
$$\mathbb{E}_{\pi_n} [   f(X_t)\mathbf{I}_{\Delta^c}(X_t)] = \mathbb{E}_{\pi_n}[ h(X_{t-1},A_{t-1})],$$
which  is the expectation of a continuous function of $(X_{t-1},A_{t-1})$ under $\mathbb{E}_{\pi_n}$. Since
$\mathbb{P}_{\pi_n}\stackrel{w}{\rightarrow}\mathbb{P}_{\pi^*}$, we get that
$$\lim_{n\rightarrow\infty} \mathbb{E}_{\pi_n} [ f(X_t)\mathbf{I}_{\Delta^c}(X_t)]=\mathbb{E}_{\pi^*}[ h(X_{t-1},A_{t-1})]
=\mathbb{E}_{\pi^*}  [f(X_t)\mathbf{I}_{\Delta^c}(X_t)].$$ The above limit trivially holds also for $t=0$.
Therefore, taking the limit as $n\rightarrow\infty$ in the rightmost expression of \eqref{eq-cond-UA-proof}
$$
\int_\XX f\mathrm{d}\mu^\XX -2\epsilon\le \mathbb{E}_{\pi^*} \Big[ \sum_{t=0}^k  f(X_t)\mathbf{I}_{\Delta^c}(X_t)\Big]\le \int_\XX f\mathrm{d}\mu_{\pi^*}^\XX$$
and so $\int f\mathrm{d}\mu^\XX \le \int f\mathrm{d}\mu_{\pi^*}^\XX$. 
To establish the reverse inequality, a similar reasoning can be done but, this time, the uniformly absorbing property is not needed. We give an outline of the proof. For $n\ge n_0$ and for arbitrary $k$ we have 
$$
\int_\XX f\mathrm{d}\mu^\XX +\epsilon \ge \int_\XX f\mathrm{d}\mu^\XX_{\pi_n}
\ge  \mathbb{E}_{\pi_n} \Big[ \sum_{t=0}^k f(X_t)\mathbf{I}_{\Delta^c}(X_t)\Big].$$
Now take the first limit as $n\rightarrow\infty$ and secondly as $k\rightarrow\infty$ to derive that 
$\int f\mathrm{d}\mu^\XX \ge \int f\mathrm{d}\mu_{\pi^*}^\XX$.
Finally, we conclude that $\int f\mathrm{d}\mu^\XX = \int f\mathrm{d}\mu_{\pi^*}^\XX$ for every bounded and continuous function $f:\XX\rightarrow[0,1]$. 
Equality of the measures $\mu^\XX=\mu_{\pi^*}^\XX$ follows.

Since $\mu$ satisfies the characteristic equations, we can disintegrate it (Lemma \ref{lemma-disintegration})  using some $\sigma\in\mathbf{M}$ so that $\mu=\mu^\XX\otimes\sigma=\mu_{\pi^*}^\XX\otimes\sigma$. Proceeding as in (3.3)--(3.4) in  \cite{DuPri24-ESAIM} we derive that for every $B\in\bcal{B}(\XX)$ 
$$0\le \mu_{\pi^*}^\XX(B)-\mu_\sigma^\XX(B)\le \int_\XX \mathbb{P}_{x,\sigma}\{T_\Delta=\infty\} \mu_{\pi^*}^\XX(dx),$$
where $\mathbb{P}_{x,\sigma}$ is the strategic probability measure of the Markovian policy $\sigma$ when the initial state of the system is $x$.
The latter integral vanishes precisely because $\mu_{\pi^*}^\XX$ is the marginal measure of an occupation measure (see \cite[Lemma 3.2]{DuPri24-ESAIM}). Hence, $\mu^\XX=\mu^\XX_{\pi^*}=\mu_\sigma^\XX$ and so $\mu=\mu_\sigma$, which is indeed an occupation measure in~$\bcal{O}$, as we wanted to prove. Note that we are not claiming that $\mu=\mu_{\pi^*}$: we have shown that $\mu=\mu_\sigma$ where $\sigma$ has been obtained from the disintegration of $\mu$. \hfill$\Box$
\\[10pt]\indent
There is an interesting side result which we state next.
\begin{proposition}
Suppose that  $\bfrak{M}$ is absorbing to $\Delta$ and that it satisfies the Conditions (W) and (SCB).
The control model $\bfrak{M}$ is uniformly absorbing to $\Delta$ if and only if the mapping $\mathbb{P}_\pi\mapsto \mathbb{E}_\pi[T_\Delta]$ is continuous on $\bcal{S}$.
\end{proposition}	
\textit{Proof.} The proof of the if part is contained in the proof of the first implication in Theorem \ref{th-W-SCB}. The only if part is derived as a straightforward consequence of the indentity \eqref{eq-tool-dugundji} combined with $\bfrak{M}$ being uniformly absorbing.
\hfill$\Box$

\subsection{Conditions on the inner and outer boundaries of $\Delta$}

We need to introduce some further notation. We define the \emph{inner boundary} of $\Delta$ as the set of points in $\partial\Delta$ which also belong to $\Delta$,  and the \emph{outer boundary} of $\Delta$ as the set of points in $\partial\Delta$ which do not belong to $\Delta$. In symbols,
$$\partial\Delta_{in}=\partial\Delta\cap\Delta\quad\hbox{and}\quad
\partial\Delta_{out}=\partial\Delta\cap\Delta^c.$$
We make the following observation: if a set $B\in\bcal{B}(\XX)$ is such that $\eta(B)=0$ and $Q(B|x,a)=0$ for every $(x,a)\in\KK$, then it is clear that the state process never visits $B$ under any policy; namely, for every $\pi\in\mathbf{\Pi}$ and every $t\ge0$, we have $\mathbb{P}_{\pi}\{X_t\in{B}\}=0$.

\begin{lemma}\label{lem-compare-mu}
Assume that  $\bfrak{M}$ is absorbing to $\Delta$ and that it satisfies the Conditions (W). Let $\pi^*,\tilde\pi,\{\pi_n\}_{n\in\NN}$ be policies such that 
$\mathbb{P}_{\pi_n}\stackrel{w}{\rightarrow} \mathbb{P}_{\pi*}$ and $
\mu_{\pi_n}\stackrel{w}{\rightarrow}\mu_{\tilde\pi}$.
Under either of the conditions
\begin{description}
\item[(i)] $\eta(\partial\Delta_{out})=0$ and $Q(\partial\Delta_{out}|x,a)=0$ for every $(x,a)\in\KK$, or
\item[(ii)] $\eta(\partial\Delta_{in})=0$, and $Q(\partial\Delta_{in}|x,a)=0$ for every $(x,a)\in\KK$, and $\bfrak{M}$ is uniformly absorbing to $\Delta$, 
\end{description}
we have that $\mu_{\pi^*}=\mu_{\tilde\pi}$.
\end{lemma}
\textbf{Proof.} We are going to show that under (i) we have $\mu_{\tilde\pi}-\mu_{\pi^*}\in\bcal{M}^+(\XX\times\AA)$, while under  (ii) we have $\mu_{\pi^*}-\mu_{\tilde\pi}\in\bcal{M}^+(\XX\times\AA)$.

Suppose now that 
 the conditions (i) are satisfied. Let $f\in\mathbf{C}_b^+(\XX\times\AA)$. Fix arbitrary  $\epsilon>0$ and $k\in\NN$. There exists $n_0$ such that for any $n\ge n_0$ we have
\begin{eqnarray}
\int_{\XX\times\AA} f\mathrm{d}\mu_{\tilde\pi}+\epsilon &\ge& \int_{\XX\times\AA} f\mathrm{d}\mu_{\pi_n}\nonumber \\
&=& \mathbb{E}_{\pi_n}\Big[ \sum_{t=0}^\infty f(X_t,A_t)\mathbf{I}_{\Delta^c}(X_t)\Big] \ =\ 
 \mathbb{E}_{\pi_n}\Big[ \sum_{t=0}^\infty f(X_t,A_t)\mathbf{I}_{\overline{\Delta}^c}(X_t)\Big]\label{eq-Delta-out}\\
& \ge&  \mathbb{E}_{\pi_n}\Big[ \sum_{t=0}^k f(X_t,A_t)\mathbf{I}_{\overline{\Delta}^c}(X_t)\Big],\nonumber
\end{eqnarray}
where the second equality in \eqref{eq-Delta-out} holds because the process $\{X_t\}$ has probability zero of visiting $\partial\Delta_{out}$ under any policy. 
Since $(x,a)\mapsto f(x,a)\mathbf{I}_{\overline{\Delta}^c}(x)$ is a lower semicontinuous function (it is the product of a nonnegative bounded continuous function and the indicator of an open set), we can take the $\liminf$ as $n\rightarrow\infty$ and obtain
$$
\int_{\XX\times\AA} f\mathrm{d}\mu_{\tilde\pi}+\epsilon \ge   \mathbb{E}_{\pi^*}\Big[ \sum_{t=0}^k f(X_t,A_t)\mathbf{I}_{\overline{\Delta}^c}(X_t)\Big] =
 \mathbb{E}_{\pi^*}\Big[ \sum_{t=0}^k f(X_t,A_t)\mathbf{I}_{\Delta^c}(X_t)\Big].$$
 But $k$ being arbitrary we obtain 
 $$
\int_{\XX\times\AA} f\mathrm{d}\mu_{\tilde\pi}+\epsilon \ge  \mathbb{E}_{\pi^*}\Big[ \sum_{t=0}^\infty f(X_t,A_t)\mathbf{I}_{\Delta^c}(X_t)\Big]=
\int_{\XX\times\AA} f\mathrm{d}\mu_{\pi^*},$$
and $\epsilon$ being arbitrary as well, we derive that  $\int f\mathrm{d}\mu_{\tilde\pi}\ge\int f\mathrm{d}\mu_{\pi^*}$. 
Hence,
 $\mu_{\tilde\pi}-\mu_{\pi^*}\in\bcal{M}^+(\XX\times\AA)$.

Assume now that the conditions (ii) hold. Given $f\in\mathbf{C}_b^+(\XX\times\AA)$ and $\epsilon>0$, let  $k$ be such that 
$$\sup_{n\in\NN} \mathbb{E}_{\pi_n}\Big[ \sum_{t>k} f(X_t,A_t)\mathbf{I}_{\Delta^c}(X_t)\Big]\le\epsilon$$
(here we use the hypothesis that $\bfrak{M}$ is uniformly absorbing). Also, there exists some 
$n_0$ such that $n\ge n_0$ implies
\begin{eqnarray*}
\int_{\XX\times\AA} f\mathrm{d}\mu_{\tilde\pi}-\epsilon &\le& \int_{\XX\times\AA} f\mathrm{d}\mu_{\pi_n}= \mathbb{E}_{\pi_n}\Big[ \sum_{t=0}^\infty f(X_t,A_t)\mathbf{I}_{\Delta^c}(X_t)\Big] \\
&\le& \mathbb{E}_{\pi_n}\Big[ \sum_{t=0}^k f(X_t,A_t)\mathbf{I}_{\Delta^c}(X_t)\Big] +\epsilon\\ &=&
 \mathbb{E}_{\pi_n}\Big[ \sum_{t=0}^k f(X_t,A_t)\mathbf{I}_{(\Delta^\circ)^c}(X_t)\Big]+\epsilon,
 \end{eqnarray*}
 the last equality being derived from the fact that the state process never visits $\partial\Delta_{in}$ under any policy. Now, since $(x,a)\mapsto f(x,a)\mathbf{I}_{(\Delta^\circ)^c}(x)$ is an upper semicontinuous function, we can take the $\limsup$ as $n\rightarrow\infty$ so that
\begin{eqnarray*}
\int_{\XX\times\AA} f\mathrm{d}\mu_{\tilde\pi}-2\epsilon &\le&  \mathbb{E}_{\pi^*}\Big[ \sum_{t=0}^k f(X_t,A_t)\mathbf{I}_{(\Delta^\circ)^c}(X_t)\Big]
= \mathbb{E}_{\pi^*}\Big[ \sum_{t=0}^k f(X_t,A_t)\mathbf{I}_{\Delta^c}(X_t)\Big]\\
&\le&  \mathbb{E}_{\pi^*}\Big[ \sum_{t=0}^\infty f(X_t,A_t)\mathbf{I}_{\Delta^c}(X_t)\Big]=\int_{\XX\times\AA} f\mathrm{d}\mu_{\pi^*}.
 \end{eqnarray*}
 We conclude that $\mu_{\pi^*}-\mu_{\tilde\pi}\in\bcal{M}^+(\XX\times\AA)$.
 
 Invoking Lemma \ref{lemma-1}, under either of the conditions (i) or (ii) above, we conclude that we indeed have the equality of the measures $\mu_{\pi^*}=\mu_{\tilde\pi}$.
\hfill$\Box$
\\[10pt]
\indent
Our next result summarizes the properties of a control model which does not visit the outer boundary of $\Delta$.  It is worth stressing that Theorem \ref{th-W-out} below does not need the Condition (SCB).
\begin{theorem}\label{th-W-out} 
Suppose that control model $\bfrak{M}$ is absorbing to $\Delta$ and it satisfies the Conditions~(W). Assume  also  that $\eta(\partial\Delta_{out})=0$ and $Q(\partial\Delta_{out} |x,a)=0$ for any $(x,a)\in\KK$. The following implications are satisfied. 
$$ \hbox{$\bfrak{P}$ is continuous} \ \Leftrightarrow\    \hbox{$\bcal{O}$ is $w$-compact}   \ \Rightarrow\ 
     \hbox{$\bfrak{M}$ is uniformly absorbing}.$$
\end{theorem}
\textit{Proof of the first equivalence.}  Obviously, if $\bfrak{P}$ is continuous then $\bcal{O}$ is $w$-compact. Regarding the reverse implication, consider policies $\{\pi_n\}_
{n\in\NN}$ and $\pi^*$ such that $\mathbb{P}_{\pi_n}\stackrel{w}{\rightarrow}\mathbb{P}_{\pi^*}$.  Since $\bcal{O}$ is a $w$-compact metric space, in order to show that $\mu_{\pi_n}\stackrel{w}{\rightarrow}\mu_{\pi^*}$ it suffices to show that $\mu_{\pi^*}$ is indeed the limit of any converging subsequence of $\{\mu_{\pi_n}\}$.  Such a result holds as a direct consequence of Lemma \ref{lem-compare-mu}(i).
\\[5pt]
\textit{Proof of the implication.}
Assume now that $\bcal{O}$ is $w$-compact.
We will prove that $\bfrak{M}$ is uniformly absorbing by contradiction. If the control model $\bfrak{M}$ is not  uniformly absorbing, there exist $\delta>0$, a sequence $\{\pi_n\}$ in $\mathbf{\Pi}$, and a sequence $\{t_n\}\uparrow\infty$ in $\NN$ such that
 \begin{equation}\label{eq-onlyif-1}
 \sum_{t>t_n} \mathbb{P}_{\pi_n} \{ T_{{\Delta}}>t\}\ge\delta\quad\hbox{for every $n\in\NN$}.
 \end{equation}
By $w$-compactness of $\bcal{S}$, there exist a subsequence of $\{\pi_n\}$ (which, without loss of generality, will be denoted again by~$\{\pi_n\}$) and a policy $\pi^*\in\mathbf{\Pi}$ such that 
$\mathbb{P}_{\pi_n}\stackrel{w}{\rightarrow}\mathbb{P}_{\pi^*}$. But $\bcal{O}$ being $w$-compact, by the first equivalence in this theorem, we also have $ \mu_{\pi_n}\stackrel{w}{\rightarrow} \mu_{\pi^*}$.
 We are going to see that this is not compatible with \eqref{eq-onlyif-1}.

Observe now that for any policy $\pi\in\mathbf{\Pi}$ we have $ \mathbf{I}_{\Delta^c}(X_t)=\mathbf{I}_{\overline{\Delta}^c}(X_t)$ with $\mathbb{P}_{\pi}$-probability one because the state process never visits $\partial\Delta_{out}$ under any policy. Therefore,
\begin{equation}\label{eq-useful}
\mu_{\pi}(\XX\times\AA)=\mathbb{E}_{\pi}\Big[\sum_{t=0}^\infty \mathbf{I}_{\Delta^c}(X_t)\Big]=
\mathbb{E}_{\pi}\Big[\sum_{t=0}^\infty \mathbf{I}_{\overline{\Delta}^c}(X_t)\Big]=
\sum_{t=0}^\infty \mathbb{P}_{\pi}\{T_{{\Delta}}>t\}.
\end{equation}
Taking $\delta>0$ from \eqref{eq-onlyif-1} and using \eqref{eq-useful} for the policy $\pi^*$
we can now find some $k$ such that 
\begin{equation}\label{eq-will-be-used}
\mu_{\pi^*}(\XX\times\AA)-\delta/2< \mathbb{E}_{\pi^*}\Big[\sum_{t=0}^k \mathbf{I}_{\overline{\Delta}^c}(X_t)\Big].
\end{equation}
Choose now $n_0$ such that $t_{n_0}\ge k$. For every $n\ge n_0$ we have 
\begin{eqnarray}
\mu_{\pi_n}(\XX\times\AA) &=&  \mathbb{E}_{\pi_n}\Big[\sum_{t=0}^{t_n} \mathbf{I}_{\overline{\Delta}^c}(X_t)\Big]+\sum_{t>t_n} \mathbb{P}_{\pi_n}\{T_{{\Delta}}>t\} \nonumber\\
&\ge &\mathbb{E}_{\pi_n} \Big[\sum_{t=0}^{k} \mathbf{I}_{\overline{\Delta}^c}(X_t)\Big]  +\delta \label{eq-onlyif-4}
\end{eqnarray}
where the last inequality comes from $t_n\ge k$ and \eqref{eq-onlyif-1}. Recalling that $\mu_{\pi_n}\stackrel{w}{\rightarrow}\mu_{\pi^*}$, we also have the convergence $\mu_{\pi_n}(\XX\times\AA)\rightarrow \mu_{\pi^*}(\XX\times\AA)$, and since $(x,a)\mapsto \mathbf{I}_{\overline{\Delta}^c}(x)$ is a lower semicontinuous function we can take the $\liminf$ in \eqref{eq-onlyif-4}
and thus
$$
\mu_{\pi^*}(\XX\times\AA) \ge  \mathbb{E}_{\pi^*}\Big[\sum_{t=0}^{k} \mathbf{I}_{\overline{\Delta}^c}(X_t)\Big]+\delta.$$
Combining this inequality with \eqref{eq-will-be-used} yields 
$\mu_{\pi^*}(\XX\times\AA)\ge \mu_{\pi^*}(\XX\times\AA)+\delta/2$, which is not possible.
 The proof of the implication is now complete.
\hfill$\Box$
\\[10pt]\indent
As a direct consequence of  Theorems \ref{th-W-SCB} and \ref{th-W-out}, we derive the following result. 
\begin{corollary}\label{cor-0}
Suppose that control model $\bfrak{M}$ is absorbing to $\Delta$ and that it satisfies the Conditions~(W) and (SCB). Assume also that $\eta(\partial\Delta_{out})=0$ and $Q(\partial\Delta_{out} |x,a)=0$ for any $(x,a)\in\KK$. Under these conditions,
$$ \hbox{$\bfrak{P}$ is continuous} \ \Leftrightarrow\       \hbox{$\bfrak{M}$ is uniformly absorbing}  \ \Leftrightarrow\ 
           \hbox{$\bcal{O}$ is $w$-compact}.$$
\end{corollary}

Now we focus on a control model which is strongly continuous at the boundary of $\Delta$ and that does not visit the inner boundary of $\Delta$.
\begin{theorem}\label{th-W-SCB-in} 
Suppose that control model $\bfrak{M}$ is absorbing to $\Delta$ and that it satisfies the Conditions~(W) and (SCB). Assume also that $\eta(\partial\Delta_{in})=0$ and $Q(\partial\Delta_{in} |x,a)=0$ for any $(x,a)\in\KK$. The following implications hold. 
$$ \hbox{$\bfrak{P}$ is continuous} \ \Leftrightarrow\       \hbox{$\bfrak{M}$ is uniformly absorbing}  \ \Rightarrow\ 
           \hbox{$\bcal{O}$ is $w$-compact}.$$
\end{theorem}
\textit{Proof.} Since we are under the conditions of Theorem \ref{th-W-SCB}, we only need to prove that when $\bfrak{M}$ is uniformly absorbing then $\bfrak{P}$ is continuous. This result, however, can be proved in the same way as the first equivalence in Theorem \ref{th-W-out} recalling that we know that $\bcal{O}$ is $w$-compact and by using, this time, item (ii) in Lemma \ref{lem-compare-mu}.\hfill$\Box$
 
\subsection{Discussion}
We have presented three results: Theorems \ref{th-W-SCB}, \ref{th-W-out}, and \ref{th-W-SCB-in}, which yield various implications between the three main issues arising in an absorbing control model $\bfrak{M}$: weak compactness of occupation measures, continuity of the natural projection mapping, and the uniformly absorbing property.  By imposing suitable topological properties on $\Delta$ or some additional conditions on the dynamic system, we can obtain several particular cases. Let us just mention the most relevant ones.
\\[5pt]\noindent
\textit{The process does not visit the boundary of $\Delta$.} If $\eta(\partial\Delta)=0$ and $Q(\partial\Delta|x,a)=0$ for every $(x,a)\in\KK$ then the process does not visit the boundary of $\Delta$ (neither, obviously, its inner nor its outer boundary) and, besides, the Condition (SCB) is also trivially satisfied because $Q(B|x,a)$ is  null  (hence, continuous) for any  measurable $B\subseteq\partial\Delta$.

\begin{corollary}\label{cor-1}
Suppose that control model $\bfrak{M}$ is absorbing to $\Delta$ and that it satisfies the Condition~(W). Assume also that $\eta(\partial\Delta)=0$ and $Q(\partial\Delta |x,a)=0$ for each $(x,a)\in\KK$. Under these conditions,
$$ \hbox{$\bfrak{P}$ is continuous} \ \Leftrightarrow\       \hbox{$\bfrak{M}$ is uniformly absorbing}  \ \Leftrightarrow\ 
           \hbox{$\bcal{O}$ is $w$-compact}.$$
\end{corollary}

The condition that $\eta(\partial\Delta)=0$ and $Q(\partial\Delta |x,a)=0$ for any $(x,a)\in\KK$ is satisfied, for instance, when $\Delta$ is an isolated set (meaning that it is at the same time an open set and a closed set) since in this case $\partial\Delta$ is the empty set.
It is  also satisfied when the MDP model is absolutely continuous with respect to some probability measure $\lambda\in\bcal{P}(\XX)$ ---by this, we mean that $\eta\ll\lambda$ and that there exists a density function $q:\XX\times\KK\rightarrow\RR^+$ such that $Q(\mathrm{d}y|x,a)=q(y|x,a)\lambda(\mathrm{d}y)$--- and the boundary of the absorbing set is such that $\lambda(\partial\Delta)=0$.
\\[5pt]
\noindent
\textit{The set $\Delta$ is closed.} If $\Delta$ is closed then $\partial\Delta_{out}$ is empty.
As a direct consequence of Theorem \ref{th-W-out} and Corollary \ref{cor-0} we have the following.

\begin{corollary}\label{cor-2}
Suppose that control model $\bfrak{M}$ is absorbing to a closed set $\Delta$ and that it satisfies the Conditions~(W). The following implications hold.
$$ \hbox{$\bfrak{P}$ is continuous} \ \Leftrightarrow\    \hbox{$\bcal{O}$ is $w$-compact}   \ \Rightarrow\ 
     \hbox{$\bfrak{M}$ is uniformly absorbing}.$$
If, in addition, the Condition (SCB) is satisfied, then all statements are equivalent.
\end{corollary}
\noindent\textit{The set $\Delta$ is open.} In this case we have that  $\partial\Delta_{in}$ is empty.
From Theorem \ref{th-W-SCB-in} we derive the following.

\begin{corollary}\label{cor-3}
Suppose that control model $\bfrak{M}$ is absorbing to an open set $\Delta$ and that it satisfies the Conditions~(W) and (SCB). The following implications hold. 
$$ \hbox{$\bfrak{P}$ is continuous} \ \Leftrightarrow\       \hbox{$\bfrak{M}$ is uniformly absorbing}  \ \Rightarrow\ 
           \hbox{$\bcal{O}$ is $w$-compact}.$$
\end{corollary}

\begin{remark}
\label{rem}
In the paper \cite{piunovskiy2023continuity} the authors study a control model $\bfrak{M}$ that is absorbing to an isolated state under the Conditions (W) or (S) ---in connection with this, see Remark \ref{rem-2} below. In \cite[Theorem 2]{piunovskiy2023continuity} it is proved that $\bfrak{P}$ is continuous if and only if $\bfrak{M}$ is uniformly absorbing, in which case, therefore, $\bcal{O}$ is $w$-compact. The authors mention as an interesting problem the following implication: if $\bcal{O}$ is $w$-compact then $\bfrak{M}$ is uniformly absorbing.

In this section, under the Conditions (W), we have given a positive answer to this conjecture in the case of an isolated $\Delta$ (cf. Corollary \ref{cor-1}). Furthermore, we have provided conditions which are weaker than $\Delta$ being isolated under which continuity of $\bfrak{P}$, $w$-compactness of $\bcal{O}$, and $\bfrak{M}$ being uniformly absorbing are all three equivalent. Namely, these conditions are that the transition kernel is strongly continuous at the boundary of $\Delta$ and that $\partial\Delta_{out}$ is never visited  (cf. Corollary \ref{cor-2}), both of which are trivially satisfied when $\Delta$ is an isolated set. 

This shows that the hypothesis that $\Delta$ is an isolated set imposes, in fact, a non-neglectable loss of generality. Allowing for a not necessarily isolated $\Delta$ and considering various properties of the transition kernel or topological properties of $\Delta$ yields a rich framework of relations between the three main statements: continuity of $\bfrak{P}$, $w$-compactness of $\bcal{O}$, and $\bfrak{M}$ being uniformly absorbing, as shown in our results in this section. Interestingly, the conditions needed to derive the various implications are closely related to the behavior of the state process at the boundary of $\Delta$ (cf. the Conditions (SCB) or the hypotheses on $\partial\Delta_{in}$ and $\partial\Delta_{out}$).
\end{remark}

\noindent
\textit{Examples.}
Our next example exhibits a weakly continuous 
control model $\bfrak{M}$ which is uniformly absorbing to a closed set $\Delta$, but whose set of occupation measures is not $w$-compact. That is, it shows a counterexample for the missing  implication in Corollary \ref{cor-2}. In view of Corollary \ref{cor-2} again, this can happen only when the Condition (SCB)  does not hold.

\begin{example}
\label{ex-1}
Let $\XX=\{-1,0\}\cup\{1/n\}_{n\ge1}$ and $\AA=\{0\}\cup\{1/n\}_{n\ge1}$. 
The spaces $\XX$ and $\AA$ (both subsets of $\RR$) are endowed with the usual topology. 
The initial state is $-1$ and $\AA(-1)=\{0\}\cup\{1/n\}_{n\ge1}$.
At the remaining states, there is just one action: $\AA(x)=\{0\}$ for $x\neq-1$. The transitions are:
$$Q(\cdot|-1,a)=\delta_a(\cdot)\quad\hbox{and}\quad Q(\cdot|x,0)=\delta_0(\cdot)\ \hbox{for $x\neq-1$}.$$
In words, at the initial state, when the controller chooses the action either $1/n$ or $0$, a transition is made to $1/n$ or $0$, respectively. Then, from that state, a transition to $0$ is made.

Let us consider the absorbing state $\Delta=\{0\}$, which is closed  but not isolated in $\XX$. Obviously, $\partial\Delta=\{0\}$ and $\partial\Delta_{out}=\emptyset$.
It is clear that the control model $\bfrak{M}$ model is uniformly absorbing  to~$\Delta$ since $T_{\Delta}\le2$  for any policy.
Conditions $(\mathrm{W})$ are satisfied: indeed, given some continuous $u\in\mathbf{C}_b(\XX)$  we have that
$$\int u(y)Q(\mathrm{d}y|x,a)=\begin{cases}
u(a)\ \hbox{for $x=-1$ and $a\in\AA(-1)$} \\
u(0)\ \hbox{for $x\neq-1$ and $a=0$}\end{cases}$$
is continuous on $\KK$, while upper semicontinuity and weak measurability of $x\mapsto\AA(x)$ also holds.

The occupation measure $\mu_{1/n}$ of the policy that chooses action $1/n$ at the state $-1$  gives mass~$1$ to $(-1,1/n)$ and mass $1$ to $(1/n,0)$. It is clear that $\mu_{1/n}$ converges weakly to a measure giving mass~$1$ to both $(-1,0)$ and $(0,0)$, which cannot be an occupation measure.  Hence, the set of occupation measures is not $w$-closed and so, it is neither $w$-compact. 
The reason why we cannot use Corollary~\ref{cor-0} is that  $Q(\partial\Delta |-1,0)=1$
and $Q(\partial\Delta |-1,1/n)=0$ and so, the transition kernel~$Q$ is not strongly continuous at the boundary of $\Delta$.
\end{example}

The following example proposes a model that satisfies the conditions of Corollary \ref{cor-0}, and so continuity of $\bfrak{P}$, $w$-compactness of $\bcal{O}$, and $\bfrak{M}$ being uniformly absorbing are all three equivalent. In this example, however, the set $\Delta$ is not isolated and so the results in \cite{piunovskiy2023continuity} cannot be invoked. 

\begin{example}
\label{ex-3}
Let $\XX=[-2,1]\cup\{2\}$ and $\AA= [0,1]$. 
The spaces $\XX$ and $\AA$ are endowed with the usual topology of $\RR$. 
Let $\Delta=[-1,0)$ and
$$\AA(x)=\begin{cases}
\{0\} & \text{ for } x\in [-2,1], \\
[0,1] & \text{ for } x=2.
\end{cases}$$
The initial distribution is given by $\eta=\delta_{\{2\}}$.
The transitions are defined by
$$Q(dy | x,a)=\begin{cases}
\nu_{\Delta}(dy) & \text{ for } (x,a)\in ](2,-1) \times\{0\}, \\
(1+x)\gamma_{\Delta}(dy)-x\nu_{\Delta}(dy) & \text{ for } (x,a)\in [-1,0) \times\{0\}, \\
\gamma_{\Delta}(dy) & \text{ for } (x,a)\in [0,1] \times\{0\}, \\
\rho_{a}(dy) & \text{ for } (x,a)\in \{2\} \times [0,1].
\end{cases}$$
where $\nu_{\Delta}$ and $\gamma_{\Delta}$ are arbitrary probability measures concentred on $\Delta$ and $\rho_{a}$ is a probability measure on $[0,1]\cup[-2,-)[$ satisfying satisfying
$\rho_{a}(\{0\})=0$ for any $a\in[0,1]$ and $\rho_{a_{n}} \stackrel{w}{\longrightarrow}\rho_{a_{*}}$ for any sequence $\{a_{n}\}_{n\in\NN}$ in $[0,1]$ converging to $a_{*}$.
It is easy to check that this model is absorbing to $\Delta$. Moreover, we have $\mathbb{P}_{\pi}\{T_{\Delta} = 2 \}=1$  for any policy $\pi$, yielding that the model is uniformly absorbing.
Conditions $(\mathrm{W})$ are satisfied. Indeed, observe that $\AA$ is compact and the multifunction $x\mapsto\AA(x)$ has  closed graph, hence it is upper semicontinuous. The action sets multifunction is also weakly measurable.
Given some continuous and bounded function $u:\XX\rightarrow\RR$ we have for any $(x,a)\in\KK$
\begin{align*}
\int_{\XX} u(y) Q(dy|x,a) & = \int_{\Delta} u(y) \nu_{\Delta}(dy) \mathbf{I}_{[-2,-1) \times\{0\}}(x,a) \\
& \quad + \Big[ (1+x)\int_{\Delta} u(y) \gamma_{\Delta}(dy) -x\int_{\Delta} u(y) \nu_{\Delta}(dy) \Big] \mathbf{I}_{[-1,0) \times\{0\}}(x,a) \\
& \quad  + \int_{\Delta} u(y) \gamma_{\Delta}(dy) \mathbf{I}_{[0,1] \times\{0\}}(x,a) + \int_{[0,1]} u(y) \rho_{a}(dy) \mathbf{I}_{\{2\} \times [0,1]}(x,a),
\end{align*}
which  is continuous in $(x,a)\in\KK$.

The occupation measure $\mu_{a}$ of the policy that chooses action $a\in[0,1]$ at the state $2$  is given by
$\mu_{a}=\delta_{(2,a)}+\rho_{a}\otimes\delta_{\{0\}}$. It is clear that the set of occupation measures $\bcal{O}=\{\mu_{a}\}_{a\in[0,1]}$ is $w$-compact
in $\bcal{M}^+(\XX\times\AA)$. The absorbing set $\Delta$ is not open nor closed and so it is not isolated. Observe that $\partial\Delta=\{-1,0\}$ with
$\partial \Delta_{out}=\{0\}$ and so, $\eta(\partial \Delta_{out})=0$ and $Q(\partial \Delta_{out} |x,a)= \rho_{a}(\{0\})=0$ for $(x,a)\in\KK$. Moreover, $\partial\Delta_{in}=\{-1\}$ and
\begin{align*}
Q(\partial\Delta_{in} |x,a) & =\nu_{\Delta}(\{-1\})\mathbf{I}_{[-2,-1]}(x)+ \big[(1+x)\gamma_{\Delta}(\{-1\})-x\nu_{\Delta}(\{-1\}) \big]\mathbf{I}_{[-1,0]}(x)\\
& \quad +\gamma_{\Delta}(\{-1\})\mathbf{I}_{[0,1]}(x)
\end{align*}
for $(x,a)\in\KK$, which is a continuous function.   This shows that the (SCB) condition is satisfied.
Consequently, the assumptions of Corollary \ref{cor-0}  hold.
\end{example}

\section{Strongly continuous absorbing control models}\label{sec-4}

In this section we assume that the control model $\bfrak{M}$ is absorbing to $\Delta$ and that it satisfies the strong continuity Conditions (S). We recall (Lemma \ref{lem-balder}) that, under these conditions, the $ws^\infty$ and the weak topologies on $\bcal{S}$ are the same, and $\bcal{S}$ is a compact metric space with these topologies.

We will need the following result. Although the proof of Lemma \ref{lem-compare-mu-S} below is   closely related to that of Lemma \ref{lem-compare-mu}, for clarity of exposition we prefer to make a detailed proof.

\begin{lemma}\label{lem-compare-mu-S}
Let $\bfrak{M}$ be absorbing to $\Delta$ and satisfy the Conditions (S). If $\pi^*,\tilde\pi,\{\pi_n\}_{n\in\NN}$ are policies such that 
$\mathbb{P}_{\pi_n}\stackrel{w}{\rightarrow} \mathbb{P}_{\pi*}$ and 
$\mu_{\pi_n}\stackrel{w}{\rightarrow}\mu_{\tilde\pi}$,
then $\mu_{\pi^*}=\mu_{\tilde\pi}$.
\end{lemma}
\textit{Proof.}  Let $f\in\mathbf{C}_b^+(\XX\times\AA)$, and fix  $\epsilon>0$ and $k\in\NN$. There is some $n_0$ such that for every $n\ge n_0$ 
\begin{eqnarray*}
\int_{\XX\times\AA} f\mathrm{d}\mu_{\tilde\pi}+\epsilon &\ge& \int_{\XX\times\AA} f\mathrm{d}\mu_{\pi_n}\nonumber \\
&=& \mathbb{E}_{\pi_n}\Big[ \sum_{t=0}^\infty f(X_t,A_t)\mathbf{I}_{\Delta^c}(X_t)\Big] \ge
 \mathbb{E}_{\pi_n}\Big[ \sum_{t=0}^k f(X_t,A_t)\mathbf{I}_{\Delta^c}(X_t)\Big].
\end{eqnarray*}
The function $(x,a)\mapsto f(x,a)\mathbf{I}_{\Delta^c}(x)$ is a bounded Carath\'eodory function and so by taking the limit as $n\rightarrow\infty$ $$
\int_{\XX\times\AA} f\mathrm{d}\mu_{\tilde\pi}+\epsilon \ge   
 \mathbb{E}_{\pi^*}\Big[ \sum_{t=0}^k f(X_t,A_t)\mathbf{I}_{\Delta^c}(X_t)\Big].$$
 Here we use the fact that $\mathbb{P}_{\pi_n}\stackrel{w}{\rightarrow} \mathbb{P}_{\pi*}$ implies $\mathbb{P}_{\pi_n}\stackrel{ws^\infty}{\longrightarrow} \mathbb{P}_{\pi*}$ since both topologies coincide for a strongly continuous control model.
Since  $k$ and $\epsilon$ were arbitrary we conclude that 
 $$
\int_{\XX\times\AA} f\mathrm{d}\mu_{\tilde\pi} \ge  \mathbb{E}_{\pi^*}\Big[ \sum_{t=0}^\infty f(X_t,A_t)\mathbf{I}_{\Delta^c}(X_t)\Big]=
\int_{\XX\times\AA} f\mathrm{d}\mu_{\pi^*}.$$
This shows that 
 $\mu_{\tilde\pi}-\mu_{\pi^*}\in\bcal{M}^+(\XX\times\AA)$ and, by Lemma \ref{lemma-1}, we necessarily have $\mu_{\tilde\pi}=\mu_{\pi^*}$.
\hfill$\Box$
\\[10pt]\indent
We state our main result in this section.

\begin{theorem}\label{th-S-main} 
Suppose that the control model $\bfrak{M}$ is absorbing to $\Delta$ and that the Conditions (S) are satisfied.  We have the following equivalences.
$$ \hbox{$\bcal{O}$ is $w$-compact} \ \Leftrightarrow\       \hbox{$\bfrak{M}$ is uniformly absorbing} \ \Leftrightarrow\  \hbox{$\bfrak{P}$ is continuous} $$
\end{theorem}
\textit{Proof.} 
\textit{If $\bcal{O}$ is $w$-compact then $\bfrak{M}$ is uniformly absorbing.} We make the proof by contradiction. Hence, 
there exist $\delta>0$, and sequences $\{\pi_n\}$ in $\mathbf{\Pi}$ and $\{t_n\}\uparrow\infty$ in $\NN$ with
 \begin{equation}\label{eq-onlyif-1-bis}
 \sum_{t>t_n} \mathbb{P}_{\pi_n} \{ T_\Delta>t\}\ge\delta\quad\hbox{for every $n\in\NN$}.
 \end{equation}
 By $w$-compactness of $\bcal{S}$ and $\bcal{O}$, there is a subsequence of $\{\pi_n\}$ still denoted by $\{\pi_n\}$ and there are policies $\pi^*,\tilde\pi\in\mathbf{\Pi}$ such that  
 $$\mathbb{P}_{\pi_n}\stackrel{w}\rightarrow\mathbb{P}_{\pi^*}\quad\hbox{and}\quad \mu_{\pi_n}\stackrel{w}{\rightarrow} \mu_{\tilde\pi}$$
 with, as a consequence of Lemma \ref{lem-compare-mu-S}, $\mu_{\pi^*}=\mu_{\tilde\pi}$.
With $\delta$ as in 
  \eqref{eq-onlyif-1-bis}, there exists some $k$ such that
\begin{equation}\label{eq-s-to-be-used}
\mu_{\pi^*}(\XX\times\AA)-\delta/2< \mathbb{E}_{\pi^*}\Big[\sum_{t=0}^k \mathbf{I}_{{\Delta}^c}(X_t)\Big].
\end{equation}
Let  $n_0$ be  such that $t_{n_0}\ge k$. Then,  $n\ge n_0$ implies
\[
\mu_{\pi_n}(\XX\times\AA) =  \mathbb{E}_{\pi_n}\Big[\sum_{t=0}^{t_n} \mathbf{I}_{{\Delta}^c}(X_t)\Big]+\sum_{t>t_n} \mathbb{P}_{\pi_n}\{T_{{\Delta}}>t\} 
\ge \mathbb{E}_{\pi_n} \Big[\sum_{t=0}^{k} \mathbf{I}_{{\Delta}^c}(X_t)\Big]  +\delta. \]
On one hand we  have  $\mu_{\pi_n}(\XX\times\AA)\rightarrow \mu_{\pi^*}(\XX\times\AA)$, and on the other hand   $(x,a)\mapsto \mathbf{I}_{{\Delta}^c}(x)$ is a Carath\'eodory function. So, by taking the limit and recalling that $\mathbb{P}_{\pi_n}$ converges to $\mathbb{P}_{\pi*}$ in the $ws^\infty$ topology,  $$
\mu_{\pi^*}(\XX\times\AA) \ge  \mathbb{E}_{\pi^*}\Big[\sum_{t=0}^{k} \mathbf{I}_{{\Delta}^c}(X_t)\Big]+\delta.$$
Combined with \eqref{eq-s-to-be-used}, this yields a contradiction, and the proof of this implication is complete. 
\\[5pt]\noindent
\textit{If $\bfrak{M}$ is uniformly absorbing then $\bfrak{P}$ is continuous.} Let $\{\pi_n\}_{n\in\NN}$ and $\pi^*$ be policies in $\mathbf{\Pi}$ such that 
$\mathbb{P}_{\pi_n}\stackrel{w}{\rightarrow} \mathbb{P}_{\pi*}$ and, therefore, also $\mathbb{P}_{\pi_n}\stackrel{ws^\infty}{\longrightarrow} \mathbb{P}_{\pi*}$. We know from Proposition \ref{prop-ua-rel-compact}  that $\bcal{O}$ is relatively $w$-compact. Hence, in order to prove that $\mu_{\pi_n}\stackrel{w}{\rightarrow}\mu_{\pi^*}$, and since $\bcal{M}^+(\XX\times\AA)$ is metrizable for the weak topology, it suffices to show that the limit of any convergent subsequence of $\{\mu_{\pi_n}\}$ is precisely $\mu_{\pi^*}$. To simplify the notation, denote also by $\{\mu_{\pi_n}\}$ the subsequence which converges to some $\mu\in\bcal{M}^+(\XX\times\AA)$.

Choose now an arbitrary  $f\in\mathbf{C}_b^+(\XX\times\AA)$. For arbitrary $k$ and $n$
\begin{eqnarray*}
\int_{\XX\times\AA} f\mathrm{d}\mu_{\pi_n}& \ge& 
 \mathbb{E}_{\pi_n}\Big[ \sum_{t=0}^k f(X_t,A_t)\mathbf{I}_{\Delta^c}(X_t)\Big].
\end{eqnarray*}
Taking the limit as $n\rightarrow\infty$ yields
\begin{eqnarray*}
\int_{\XX\times\AA} f\mathrm{d}\mu& \ge& 
 \mathbb{E}_{\pi^*}\Big[ \sum_{t=0}^k f(X_t,A_t)\mathbf{I}_{\Delta^c}(X_t)\Big]
\end{eqnarray*}
and thus $\int f\mathrm{d}\mu\ge \int f\mathrm{d}\mu_{\pi^*}$. Now, since $\bfrak{M}$ is uniformly absorbing, for any $\epsilon$ there exists some $k$ such that, for every $n\in\NN$,
\begin{eqnarray*}
\int_{\XX\times\AA} f\mathrm{d}\mu_{\pi_n}& \le& 
 \mathbb{E}_{\pi_n}\Big[ \sum_{t=0}^k f(X_t,A_t)\mathbf{I}_{\Delta^c}(X_t)\Big]+\epsilon.
\end{eqnarray*}
Taking again the limit as $n\rightarrow$ gives
\begin{eqnarray*}
\int_{\XX\times\AA} f\mathrm{d}\mu& \le& 
 \mathbb{E}_{\pi^*}\Big[ \sum_{t=0}^k f(X_t,A_t)\mathbf{I}_{\Delta^c}(X_t)\Big]+\epsilon\ \le\  \int_{\XX\times\AA} f\mathrm{d}\mu_{\pi^*}+\epsilon.
\end{eqnarray*}
Therefore, $\int f\mathrm{d}\mu\le \int f\mathrm{d}\mu_{\pi^*}$. This shows that we have $\mu=\mu_{\pi^*}$, which completes the proof.
\\[5pt]\noindent
\textit{If  $\bfrak{P}$ is continuous then $\bcal{O}$ is $w$-compact.} This is the obvious implication.
\hfill$\Box$
\\[10pt]\indent
It is worth stressing that we have obtained here, under the Conditions (S), the equivalence of all the three main properties  regardless of the topological properties of $\Delta$.  

Recall now  that the  model in Example \ref{ex-1}  is uniformly absorbing but $\bcal{O}$ is not $w$-compact. It does not verify, however, the Conditions (S) since $Q(\{0\}| -1,1/n)=0$ does not converge to $Q(\{0\}|-1,0)=1$.
 So, we cannot use Theorem \ref{th-S-main}.

\begin{remark}
\label{rem-2}
In connection with Remark \ref{rem}, we recall that 
the paper \cite{piunovskiy2023continuity} also considers a control model $\bfrak{M}$ under the Conditions (S) which is absorbing to an isolated set $\Delta$. In Theorem \ref{th-S-main} we have again given a positive answer  to the conjecture that when $\bcal{O}$ is $w$-compact then $\bfrak{M}$ is uniformly absorbing. 
 What is most interesting, however, is that we have proved equivalence of  
continuity of~$\bfrak{P}$, $w$-compactness of $\bcal{O}$, and $\bfrak{M}$ being uniformly absorbing without any  topological assumption  on  $\Delta$: indeed, we just require that $\Delta$ is a measurable subset of the state space $\XX$. 
\end{remark}

\paragraph{The $ws$-topology.}

In this paper we have focused on the weak topology on the set of finite measures $\bcal{M}^+(\XX\times\AA)$. However, the paper \cite{DuPri24-ESAIM} studied  absorbing control model using the $ws$-topology for the set $\bcal{O}$ of occupation measures. With the results in the present paper, we can obtain some further results linking both topologies.

We recall that the $ws$-topology on $\bcal{M}^+(\XX\times\AA)$ is the coarsest topology for which the mappings $\mu\mapsto \int f\mathrm{d}\mu$ are continuous for every bounded Carath\'eodory function $f:\XX\times\AA\rightarrow\RR$, i.e., $f$ is measurable and $f(x,\cdot)$ is continuous on $\AA$ for each fixed $x\in\XX$. It is clear that the $ws$-topology is finer than the weak topology.

The result in Proposition \ref{prop-ws} below was proved in \cite{DuPri24-ESAIM}. Here we provide an alternative proof based on our results in the present paper and on  Sch\"{a}l's criterion for relative compactness in the $ws$-topology, which we state next.

\begin{lemma}\label{lem-schal-ws} Let $\bfrak{M}$ be absorbing to $\Delta$. The set $\bcal{O}$ is relatively $ws$-compact if and only if  for every sequence $\{f_n\}$ of bounded Carath\'eodory functions on $\XX\times\AA$ such that $f_n\downarrow0$ we have
$$\lim_{n \rightarrow\infty}\sup_{\mu\in\bcal{O}} \int_{\XX\times\AA} f_n \mathrm{d}\mu =0.$$
\end{lemma}
\textit{Proof.} This result is proved in \cite[Theorem 3.10]{schal75} for sets of probability measures in $\bcal{P}(\XX\times\AA)$. However, assuming that $\eta(\Delta^c)>0$   (otherwise, all occupation measures would be null) and recalling that $\bcal{O}$ is bounded, it follows that there are constants such that $0<m\le \mu(\XX\times\AA)\le M$ for every $\mu\in\bcal{O}$, and the result in \cite{schal75} can be easily adapted to obtain the lemma.
\hfill$\Box$
%

\begin{proposition}\label{prop-ws} Let the control model $\bfrak{M}$ be absorbing to $\Delta$ and satisfy the Conditions (S). The set of occupation measures $\bcal{O}$ is $ws$-compact if and only if $\bfrak{M}$ is uniformly absorbing.
\end{proposition}
\textit{Proof.} It is clear that when $\bcal{O}$ is $ws$-compact, then $\bcal{O}$ is $w$-compact. By Theorem \ref{th-S-main}, if $\bcal{O}$ is $w$-compact then $\bfrak{M}$ is uniformly absorbing. Hence, it remains to show that when $\bfrak{M}$ is uniformly absorbing, we have that $\bcal{O}$ is $ws$-compact. 
To establish that $\bcal{O}$ is relatively $ws$-compact we can use the criterion in Lemma \ref{lem-schal-ws}: indeed, just repeat the proof of Lemma \ref{prop-ua-rel-compact}  taking bounded Carath\'eodory functions $f_n\downarrow0$ and using the $ws^\infty$-topology on $\bcal{S}$, for which $\bcal{S}$ is compact. 
Hence, it remains to show that $\bcal{O}$ is $ws$-closed. If $\{\mu_\alpha\}$ is a net in $\bcal{O}$ converging in the $ws$-topology to some $\mu\in\bcal{M}^+(\XX\times\AA)$, convergence in the $w$-topology also takes place. But $\bcal{O}$ being $w$-compact ($\bfrak{M}$ is uniformly absorbing and we use Theorem \ref{th-S-main}), we necessarily have $\mu\in\bcal{O}$.
\hfill$\Box$
\\[10pt]\indent
We derive the following result.
\begin{corollary}\label{cor-ws}
Let $\bfrak{M}$ be uniformly absorbing to $\Delta$ and satisfy the Conditions (S). The weak  topology and the  $ws$-topology on $\bcal{O}$ coincide. 
\end{corollary}
\textit{Proof.} 
From \cite[Corollary 2.2 and Proposition 2.3]{balder01}, we derive that the $ws$-topology is metrizable on~$\bcal{O}$. But $\bcal{O}$ being compact under both the weak  and the $ws$-topologies, and the $ws$-topology being finer than the weak topology, the stated result easily follows. 
\hfill$\Box$
\\[10pt] 
\indent We stress again that the results in Proposition \ref{prop-ws} and Corollary \ref{cor-ws} are obtained without any topological assumption on the set $\Delta$.

\end{document}